\documentclass[preprint,3p]{elsarticle}
\usepackage{amssymb,amscd,amsxtra,calc}
\usepackage{cmmib57}
\usepackage{multirow}
\usepackage{mathrsfs}
\usepackage{graphicx}
\usepackage{bm}
\usepackage{amsmath}
\usepackage{amsthm}
\usepackage{amssymb}
\usepackage{graphicx}
\usepackage{accents}
\usepackage{tikz} %pictures
\usepackage{tikz-cd}%implies
\usepackage{pgfplots}
\pgfplotsset{compat=1.18}
\usepgfplotslibrary{groupplots}
\definecolor{seaborngreen}{rgb}{0.3333333333333333, 0.6588235294117647, 0.40784313725490196}
\definecolor{seabornblue}{rgb}{0.2980392156862745, 0.4470588235294118, 0.6901960784313725}
\definecolor{seabornred}{rgb}{0.7686274509803922, 0.3058823529411765, 0.3215686274509804}
\usepackage{stmaryrd}
\usepackage{comment}
\usepackage{colonequals}
\usepackage{subcaption}
\usepackage[section]{placeins}
\usepackage{algorithmic}
\usepackage{algorithm}
\usepackage{graphicx}
\allowdisplaybreaks[3]
\usepackage{hyperref}
\usepackage{aliascnt}
\usepackage[capitalize,nameinlink]{cleveref}
\hypersetup{
    colorlinks,
    allcolors={green!50!black},
    urlcolor={red!50!black}
}

\usepackage{colonequals}

\theoremstyle{plain}
\newtheorem{theorem}{Theorem}[section]
\newaliascnt{lemma}{theorem}
\newtheorem{lemma}[lemma]{Lemma}
\aliascntresetthe{lemma}
\newaliascnt{corollary}{theorem}

\aliascntresetthe{corollary}
\newaliascnt{proposition}{theorem}

\aliascntresetthe{proposition}
\theoremstyle{definition}
\newaliascnt{definition}{theorem}
\newtheorem{definition}[definition]{Definition}
\aliascntresetthe{definition}
\theoremstyle{remark}
\newaliascnt{example}{theorem}
\newtheorem{example}[example]{Example}
\aliascntresetthe{example}
\newaliascnt{remark}{theorem}
\newtheorem{remark}[remark]{Remark}
\aliascntresetthe{remark}
\crefname{example}{Example}{Examples}
\Crefname{example}{Example}{Examples}
\crefname{remark}{Remark}{Remarks}
\Crefname{remark}{Remark}{Remarks}

\AddToHook{env/theorem/begin}{\crefalias{section}{theorem}}
\AddToHook{env/lemma/begin}{\crefalias{section}{lemma}}
\AddToHook{env/corollary/begin}{\crefalias{section}{corollary}}
\AddToHook{env/proposition/begin}{\crefalias{section}{proposition}}
\AddToHook{env/definition/begin}{\crefalias{section}{definition}}
\AddToHook{env/example/begin}{\crefalias{section}{example}}
\AddToHook{env/remark/begin}{\crefalias{section}{remark}}

\numberwithin{equation}{section}
\numberwithin{algorithm}{section}

\def\Proof{\noindent{\bf Proof} \quad}
\def\qed{\hfill $\Box$ \smallskip}

\renewcommand{\d}{\mathrm{d}}
\newcommand{\p}{\partial}

\renewcommand{\l}{\left}
\renewcommand{\r}{\right}

\def\R{{\mathbb R}}

\def\cM{{\mathcal M}}
\def\cC{{\mathcal C}}

\renewcommand{\S}{\mathbb{S}}

\newcommand{\cE}{\mathcal{E}}

\newcommand{\I}{\mathbb{I}}

\journal{Journal of Computational Physics}

\begin{document}
\begin{frontmatter}

\title{A minimizing-movement framework for parametric finite element approximation of anisotropic curve shortening flows}

\author[author1]{Wei Jiang}
\ead{jiangwei1007@whu.edu.cn}
\address[author1]{School of Mathematics and Statistics, Hubei Key Laboratory of Computational Science, Wuhan University, Wuhan 430072, China}

\author[author2]{Chunmei Su}
\ead{sucm@tsinghua.edu.cn}
\address[author2]{Yau Mathematical Sciences Center, Tsinghua University, Beijing 100084, China}
\cortext[cor1]{Corresponding author.}

\author[author3]{Ganghui Zhang\corref{cor1}}
\ead{ganghui.zhang@polyu.edu.hk}
\address[author3]{Department of Mathematics, The Hong Kong Polytechnic University, Hong Kong, China}

\begin{abstract}
We propose a minimizing-movement framework for  constructing
parametric finite element methods (PFEMs) for anisotropic curve shortening flows.  The cornerstone is an \emph{admissible modified energy}: a smooth, convex, and iteration-dependent surrogate that majorizes the original anisotropic energy and agrees with it at the current iterate.
For any fixed finite element degree $n\geq1$ and any positive surface energy density $\gamma\in C^2(\mathbb{S}^1)$, the unconstrained one-step problem has a unique solution for any time-step size whenever the parametric derivative of the current discrete curve does not vanish on any mesh element. It also satisfies a one-step discrete energy dissipation inequality. It therefore preserves the variational structure of the underlying $L^2$-gradient flow. The framework recovers a variational interpretation of Dziuk's isotropic scheme  and applies directly to the weak anisotropies. For strong anisotropies, a convex splitting removes the weak-anisotropy restriction on the original density. A scalar constraint further extends the framework to area-preserving flows.
Numerical experiments exhibit approximately the expected convergence for smooth isotropic and weakly anisotropic examples. They also confirm energy dissipation and area preservation. In the strongly anisotropic regime, the computations produce sharp facet-like profiles but reveal mesh-scale oscillations under the grid refinement.
\end{abstract}

\begin{keyword}
anisotropic curve shortening flow \sep parametric finite element method \sep
minimizing-movement \sep energy stability \sep convex splitting
\MSC[2020] 65M60 \sep 65M12 \sep 35K93
\end{keyword}

\end{frontmatter}

%\tableofcontents

\section{Introduction}\label{sec:intro}
Geometric flows describe the evolution of curves and surfaces driven by
intrinsic geometric quantities. They arise in materials science
\cite{garcke2013curvature,cahn1994overview}, image processing
\cite{sapiro2006geometric}, and computational geometry
\cite{clarenz2000anisotropic}.  For a family of closed curves $\Gamma(t)$ in
$\R^2$, the isotropic \emph{curve shortening flow} (CSF) is given by
$$V_n=-\kappa,$$
where $V_n$ and $\kappa$ denote the normal velocity and curvature, respectively. CSF is the $L^2$-gradient flow
of the perimeter $\int_\Gamma \d s$.

In crystalline materials, the interfacial energy depends on the interface orientation.  This dependence is important in crystal growth \cite{taylor1992overview,gurtin1993thermomechanics}, grain-boundary
migration \cite{kinderlehrer2001evolution}, and solid-state dewetting
\cite{Bao2017}. The corresponding free energy is the \emph{anisotropic} surface
energy \cite{Taylor1992,Taylor1994,hoffman1972vector,cahn1974vector}
\begin{equation}\label{eq:intro-aniso-energy}
\cE_{\gamma}(\Gamma) := \int_\Gamma \gamma(N)\,\d s,
\end{equation}
where $\gamma : \S^1 \to \R^+=(0,\infty)$ is the surface
energy density and $N$ is the outward unit normal to $\Gamma$. The
associated $L^2$-gradient flow is the \emph{anisotropic curve
shortening flow}. In this flow,
$\kappa$ is replaced by the weighted curvature $\kappa_\gamma$, which
depends on $\gamma$ and its second derivative with respect to orientation. Let $\widehat{\gamma}$ denote the $1$-homogeneous extension of $\gamma$.
The density $\gamma$ is called \emph{weakly anisotropic}
\cite{Bao2017} if
\[
\nabla^2\widehat{\gamma}(N)T\cdot T\geq0
\]
for every $N\in\S^1$ and every unit tangent vector $T\perp N$. In this regime, the anisotropic curvature operator is parabolic. If this condition fails, $\gamma$ is called \emph{strongly anisotropic}, and the operator loses parabolicity in some orientations; the unregularized flow may be ill-posed
\cite{gurtin1993thermomechanics,angenent1989multiphase}.  Strong anisotropies are also associated with missing orientations and may produce Wulff shapes with facets and corners \cite{Bao17stable}. Their analysis and numerical
approximation are therefore considerably more delicate
\cite{Bellettini2001,Giga2006,Bao2017}. A central computational challenge
is to construct accurate and energy-stable schemes that apply to both anisotropic regimes.

The parametric finite element method (PFEM) directly tracks an evolving interface and has become an important tool for simulating geometric flows.
Dziuk introduced seminal PFEMs for mean curvature flow \cite{Dziuk1991} and anisotropic curve shortening flow \cite{Dziuk1994, Dziuk1999}.
Barrett, Garcke, and N{\"u}rnberg subsequently developed
a family of variational PFEMs for several second- and fourth-order geometric flows \cite{BGN07A,BGN07B,BGN20}. These methods combine discrete energy stability with favorable tangential redistribution of mesh
points. The surveys \cite{Deckelnick2005,BGN20} provide broader accounts of these developments.  Recent studies have also advanced the analysis of PFEMs through evolving surface techniques
\cite{Kovacs-Li-Lubich2021}, artificial tangential motions \cite{Hu-Li}, monotone structures \cite{Li1}, projected distance methods \cite{Bai-Hu-Li}, and stabilized BGN-type formulations \cite{Bai2023}.

Considerable effort has been devoted to extending energy-stable PFEMs to anisotropic geometric flows. Barrett {\it et al.} treated anisotropies induced
by Riemannian metrics \cite{BGN2008c,BGN2008d}. Bao {\it et al.} incorporated a Willmore energy
regularization to simulate solid-state dewetting with
strong anisotropy~\cite{Bao17stable,Bao2017}.  Bao and Li later developed an energy-stable PFEM for planar curves governed by more general anisotropic surface energies~\cite{Li2021energy}. Their stability analysis, however, requires anisotropy-dependent structural conditions. To overcome this limitation, Bao, Jiang, and Li then introduced a symmetrized surface energy matrix and constructed a symmetrized energy-stable PFEM~\cite{Bao2023b}, which is energy stable for surface energies satisfying $$\gamma(N)=\gamma(-N).$$
A subsequent surface-energy-matrix construction relaxed this restriction~\cite{Bao2024}
to $$3\gamma(N)>\gamma(-N), \qquad \forall \, N\in\S^1.$$
These developments substantially enlarge the admissible
class of surface energies. Nevertheless, their one-step energy estimates still depend on pointwise conditions involving the anisotropy.

This limitation motivates the following question: \emph{ can one construct an energy-stable PFEM for anisotropic CSF that applies to every positive $\gamma\in C^2(\S^1)$ without imposing symmetry or a pointwise comparison between $\gamma(N)$ and $\gamma(-N)$?} The present
work answers this question affirmatively. We introduce a variational framework based on an iteration-dependent modified energy. Its principal
advantage over existing surface-energy-matrix PFEMs is the broader scope of its stability guarantee. In particular, the proposed analysis covers nonsymmetric and strongly anisotropic densities that may violate
$3\gamma(N)>\gamma(-N)$.

 Our construction is motivated by the \emph{minimizing-movement scheme} (MMS) introduced by De~Giorgi \cite{DeGiorgi1993}. Given an energy
functional $\cE$ on a metric space $(\cM,d)$ and a current state
$u^k\in\cM$, the classical MMS defines the next state by
\begin{equation}\label{eq:intro-abstract-MMS}
u^{k+1}
\in
\mathop{\arg\min}_{u\in\cM}
\left\{
\frac{1}{2\tau}d(u,u^k)^2+\cE(u)
\right\}.
\end{equation}
When a minimizer exists, comparison with $u^k$ immediately gives an energy
inequality. Existence and uniqueness, however, require additional
compactness, coercivity, or convexity assumptions.
The MMS has played an important role in several areas of applied analysis. It underlies the Wasserstein gradient-flow formulation of the
Fokker--Planck equation introduced by Jordan, Kinderlehrer, and Otto \cite{JKO}. It also provides a foundation for weak formulations of mean curvature flow \cite{ATW,Luckhaus,Chambolle2004} and motivates many
energy-stable phase-field schemes \cite{Du-Feng2020}. We refer to \cite{Ambrosio} for a comprehensive treatment.

For curvature-driven interfaces, classical minimizing-movement formulations typically represent the evolving geometry by phase fields or
sets of finite perimeter. Their distance terms are well suited to
weak-solution theory but do not directly yield the standard variational
forms used by parametric finite element discretizations. A recent parametric framework by Liu and
Zhao \cite{LiuZhao2026} treats isotropic mean curvature flow and surface
diffusion. Their method combines a metric penalty for normal displacement
with weak constraints on tangential motion. It recovers the BGN and dual
minimal-deformation-rate schemes and produces mesh-improving variants. The
present work addresses a different difficulty. We focus on the nonconvex
anisotropic perimeter and seek an energy estimate that does not require
additional conditions on $\gamma$.

To connect minimizing-movements with PFEMs, let
$\Gamma_h^k=X_h^k(\I)$ be the approximation curve at time $t_k=k\tau$. A natural time-discrete update is
\begin{equation}\label{eq:intro-MMS}
X_h^{k+1}
\in
\mathop{\arg\min}_{X\in[V_h]^2}
\left\{
\frac{1}{2\tau}
\|X-X_h^k\|_{L^2(\Gamma_h^k)}^2
+
\cE_\gamma(X)
\right\},
\end{equation}
where $V_h$ is the continuous periodic piecewise-polynomial finite element space of a fixed degree $n\geq1$ on the parameter domain $\I$.
The first term is a weighted $L^2$-penalty induced by the previous discrete parametric curve. The second term is the anisotropic energy of the curve $X(\I)$. This choice avoids a bulk distance functional and is directly compatible with the PFEM representation.

Under the elementwise non-degeneration   assumption on the previous discrete parametrization, the quadratic penalty term in \eqref{eq:intro-MMS} is coercive with respect to the full vector of finite element coefficients. The finite-dimensional objective therefore has a global
minimizer. Directly solving this problem is nevertheless difficult. Its
Euler--Lagrange equation contains the full nonlinear gradient of
$\cE_\gamma$ and produces a strongly coupled implicit system. Moreover,
the objective may be nonconvex in the strongly anisotropic regime. A
minimizer need not be unique, and the stationarity equations may have
multiple solutions. Newton-type methods may therefore be sensitive to their initial guesses.

We overcome these by introducing an
\emph{admissible modified energy}
$\cE_{\mathrm{m},\gamma}^k(X)$. This functional is smooth, convex, and
dependent on the previous iterate. It satisfies the consistency condition
\[
D\cE_{\mathrm{m},\gamma}^k(X_h^k)
=
D\cE_\gamma(X_h^k)
\]
and the majorization condition required to control the true anisotropic
energy. Replacing $\cE_\gamma$ in \eqref{eq:intro-MMS} with the admissible
modified energy gives a strictly convex time-step problem. Consequently,
the discrete update is unique and satisfies
\[
\cE_\gamma(X_h^{k+1})
+
\frac{1}{2\tau}
\|X_h^{k+1}-X_h^k\|_{L^2(\Gamma_h^k)}^2
\leq
\cE_\gamma(X_h^k).
\]

An admissible modified energy can be constructed directly in the weakly
anisotropic regime. For a strongly anisotropic density, we use the
convex-splitting decomposition \cite{Eyre1998,dai2022convex}
\[
\gamma=\gamma_c-\gamma_e.
\] We treat the convex contribution associated
with $\gamma_c$ implicitly and the concave contribution
$-\gamma_e$ explicitly. This construction applies to every positive
density $\gamma\in C^2(\S^1)$. The resulting method therefore satisfies
the discrete energy inequality at every time step without requiring
symmetry or an additional pointwise condition on $\gamma$.

The main contributions of this work are summarized as follows.
\begin{itemize}
\item
We formulate a parametric-finite-element-compatible minimizing-movement principle using a
weighted $L^2$-penalty on the previous discrete curve. This construction
bridges the gap between set-based minimizing movements and explicit PFEM representations.

\item
We introduce an admissible modified energy that yields a uniquely solvable
time-step problem and a one-step dissipation inequality for the true
anisotropic energy.

\item
We combine the modified-energy construction with convex splitting. The
resulting PFEM applies to every positive
$\gamma\in C^2(\S^1)$ without requiring
$\gamma(N)=\gamma(-N)$ or
$3\gamma(N)>\gamma(-N)$. \emph{To the best of our knowledge, this is the first fully discrete PFEM for anisotropic curve shortening flow with a one-step
energy estimate for this entire class of surface energy densities}.

\item
We extend the framework to area-preserving curve shortening flow through a
Lagrange multiplier formulation. We also show that the isotropic specialization recovers Dziuk's classical PFEM \cite{Dziuk1994} as the
Euler--Lagrange equation of one minimizing-movement step. This result gives
a new incremental variational interpretation of Dziuk's scheme.
\end{itemize}

The proposed construction is inspired by the MMS but is not same with the classical
De~Giorgi scheme for the fixed energy $\cE_\gamma$. The modified energy
depends on the previous iterate. Nevertheless, both methods compute the
next state by minimizing a time-increment penalty plus an energy term. The
first-variation consistency condition ensures that the modified problem
retains the local variational structure of the original anisotropic
energy.

The remainder of the paper is organized as follows.
\Cref{sec:mms-framework} develops the abstract minimizing-movement framework and the admissible modified-energy construction. It first
applies the framework to isotropic CSF and recovers Dziuk's classical
PFEM. It then extends the construction to weakly and strongly anisotropic
surface energies. \Cref{sec:ap-csf} treats area-preserving CSF through a
Lagrange multiplier formulation. \Cref{sec:numerics} presents numerical
experiments that examine accuracy, energy dissipation, area preservation,
solver robustness, and mesh quality. Finally, \Cref{sec:conclusion}
summarizes the main results and discusses possible extensions.

\section{A minimizing-movement framework for CSF}
\label{sec:mms-framework}

In this section, we develop a unified variational framework
inspired by De~Giorgi's MMS.  We first formulate an
abstract framework which gives unique one-step solvability and an energy estimate for any
time-step size.
We then apply the framework to three concrete
cases. The isotropic construction recovers Dziuk's classical semi-implicit scheme. The second construction treats weakly anisotropic
energies and the third combines the framework with convex splitting to treat
strongly anisotropic energies.

Throughout this section, let $\I=[0,1]$ be the periodic reference interval, with its endpoints identified, and consider the uniform partition
\[
\I=\bigcup_{j=1}^J I_j,
\qquad
I_j=[\xi_{j-1},\xi_j],
\qquad
\xi_j=jh,
\qquad
h=\frac1J.
\]
Fix an integer $n\geq1$ and define
\[
V_h^n
:=\bigl\{
v_h\in C^0(\I,\R):
v_h|_{I_j}\in\mathbb P^n(I_j),\;
1\leq j\leq J,\;
v_h(0)=v_h(1)
\bigr\},
\]
where $\mathbb P^n(I_j)$ denotes the space of real-valued polynomials of degree at most $n$ on $I_j$.
For brevity, we write $V_h:=V_h^n$ below.
Given $X_h^k\in[V_h]^2$, set $\Gamma_h^k:=X_h^k(\I)$ and assume that its elementwise derivative satisfies $|\p_\xi X_h^k|>0$ on every element.
We seek the next discrete parametrization $X_h^{k+1}\in[V_h]^2$.
For $f\in[V_h]^2$, define
\[
\|f\|_{L^2(\Gamma_h^k)}^2
:=\int_\I |f|^2\,|\p_\xi X_h^k|\,\d\xi.
\]
All variational integrals in the analysis are understood as sums of exact integrals over the elements.

Let $N_h:=\dim V_h$.
After fixing a basis of $V_h$, we identify
\[
[V_h]^2\cong\R^{2N_h}
\]
by associating each vector-valued finite element function with its full coefficient vector.
Accordingly, convexity, continuity, coercivity, and differentiability of functionals on $[V_h]^2$ are understood in the usual finite-dimensional sense under this identification.
If
$\cE:[V_h]^2\to\R$ is continuously differentiable, we denote its first
variation by
\begin{equation}\label{eq:first-variation-def}
D\cE(X)(v)
:=
\left.
\frac{\d}{\d\varepsilon}
\right|_{\varepsilon=0}
\cE(X+\varepsilon v),
\qquad
v\in[V_h]^2.
\end{equation}

\subsection{An abstract minimizing-movement framework}
\label{sec:mms-abstract}
Let $X: \I\rightarrow \R^2$ be a counterclockwise parametrization of a closed planar curve, and set
\(\Gamma(X):=X(\I)\). Writing $s$ for the arclength parameter, we define the unit tangent and the outward unit normal by
\[T=\p_sX, \qquad N=T^\perp,
\]
where $(\cdot)^\perp$ denotes clockwise rotation through $\pi/2$. With this convention, the  curvature \(\kappa\) is defined by the Frenet relation $\p_sT=-\kappa N$.

Let $\gamma\in C^2(\S^1)$ be a positive surface energy density, and let
\[
\widehat\gamma(p)
=
|p|\gamma\left(p/|p|\right),
\qquad
p\in\R^2\setminus\{0\},
\]
denote its positively one-homogeneous extension. We set $\widehat\gamma(0)=0$. Thus
\(\widehat\gamma\in C^2(\R^2\setminus\{0\})\).
The corresponding anisotropic interfacial energy is
\[
\cE_\gamma(X)
:=
\int_{\Gamma(X)}\gamma(N)\,\d s
=\int_\I
\widehat\gamma\bigl((\p_\xi X)^\perp\bigr)\,\d\xi.
\]
Modulo tangential reparametrizations, the $L^2$-gradient flow of $\cE_\gamma$ is the anisotropic CSF:
$$\p_t X=-\kappa_\gamma N,$$
where the anisotropic curvature is given by
$$\kappa_\gamma=\nabla^2\widehat\gamma(N)T\cdot T\,\kappa.$$
In the isotropic case $\gamma\equiv 1$, one has \(\widehat\gamma(p)=|p|\) and
$$ \nabla^2\widehat\gamma(N)T\cdot T=1. $$
Consequently, \(\kappa_\gamma=\kappa\), and the evolution reduces to the classical CSF.
%For convenience, we write
%\[ p(X):=(\p_\xi X)^\perp, \qquad p^k:=(\p_\xi X_h^k)^\perp.\]

\begin{definition}
The density $\gamma$ is called \emph{weakly anisotropic} if its {\emph{stiffness}} satisfies
\begin{equation}\label{eq:weak-aniso}
\nabla^2 \,\widehat\gamma(N)\,T \cdot T \geq 0
\qquad \forall\, N \in \mathbb{S}^1,\quad T\perp N.
\end{equation}
 Otherwise, $\gamma(N)$ is called \emph{strongly anisotropic}.
\end{definition}
 In angular coordinates, if
$$ N(\theta)=(\cos\theta,\sin\theta), \qquad T(\theta)=(-\sin\theta,\cos\theta), $$
then the weak anisotropy condition is equivalently
\[\nabla^2 \,\widehat\gamma(N)\,T \cdot T
=\gamma(\theta)+\gamma''(\theta)\ge 0.\]

We use the non-strict condition \eqref{eq:weak-aniso}. This definition
includes densities whose stiffness vanishes at isolated orientations, such
as the $\ell^4$-norm anisotropy considered below. Under the non-strict condition, the continuous equation may be degenerate parabolic. Uniform
parabolicity requires a strictly positive lower bound
\[
\nabla^2\widehat\gamma(N)T\cdot T\geq c_\gamma>0.
\]

\begin{example}\label[example]{Example:Ani}
The following surface energy densities are commonly used in the literature.
\begin{itemize}
\item \textbf{Case 1 (isotropic).} Let $\gamma(N) \equiv 1$.  Then
$\nabla^2\,\widehat\gamma(N)\,T \cdot T = 1$.

\item \textbf{Case 2 ($q$-fold anisotropy).}
Let
$$\gamma(N) = 1 + \beta\cos(q\theta),\qquad N = (-\sin\theta,\cos\theta)^\top \in \mathbb{S}^1,$$ where
$q\geq2$ is an integer and $0\leq\beta<1$. The restriction on $\beta$ ensures that $\gamma$ is positive. A direct calculation yields
$$\nabla^2 \,\widehat\gamma(N)\,T \cdot T
= 1 - \beta(q^2 - 1)\cos(q\theta).$$
Consequently, the density is weakly
anisotropic when $\beta \leq 1/(q^2-1)$ and strongly anisotropic when $\beta > 1/(q^2-1)$~\cite{Bao2017}.

\item \textbf{Case 3 (Riemannian-type anisotropy).}
Let $$\gamma(N) = \sum_{\ell=1}^M \sqrt{N^\top G_\ell N},$$
 where each
$G_\ell \in \R^{2\times 2}$ is symmetric positive definite. This density is weakly anisotropic in the sense of \eqref{eq:weak-aniso}; see \cite{BGN2008c,Bao2023b}.

\item \textbf{Case 4 ($\ell^r$-norm anisotropy).}
Let $$\gamma(N) = (|n_1|^r + |n_2|^r )^{1/r},\qquad N = (n_1,n_2)^\top \in \mathbb{S}^1,\qquad 2\leq r<\infty.$$
 This density is weakly anisotropic~\cite{Bao2017,Bao2023b}.
For $1<r<2$, it is not $C^2$ at orientations
on the coordinate axes. It therefore lies beyond the present analysis, although the discrete formulation can still be implemented and examined
numerically.
\end{itemize}
\end{example}

\smallskip

Let $\cE_{\gamma} : [V_h]^2 \to \R$ be a continuous energy functional. A direct minimizing-movement discretization is given by
\begin{equation}\label{eq:MMS-direct}
X_h^{k+1} \in \mathop{\arg\min}_{X \in [V_h]^2}
\biggl\{ \frac{1}{2\tau}\|X - X_h^k\|_{L^2(\Gamma_h^k)}^2
+ \cE_{\gamma}(X) \biggr\}.
\end{equation}
Direct calculation gives the first variation of $\cE_{\gamma}$:
 \begin{align*}
    D \cE_{\gamma}(X)(v)
    &=\l.\frac{\d }{\d \varepsilon}\r|_{ \varepsilon=0}\int_\I \widehat\gamma\bigl((\p_\xi X+\varepsilon \p_\xi v)^\perp\bigr)\,\d\xi=\int_\I \nabla \widehat\gamma\bigl((\p_\xi X)^\perp\bigr)\cdot (\p_\xi v)^\perp \,\d\xi.
\end{align*}
Therefore, if $\widehat\gamma$ is not differentiable at the origin, the integrand in this formula is undefined at points where
$\p_\xi X=0$. Although isolated zeros of $\p_\xi X$ do not necessarily preclude differentiability of the integrated energy, the objective in
\eqref{eq:MMS-direct} need not be continuously differentiable throughout the trial space. The direct formulation thus gives rise to a nonlinear, potentially
nonsmooth optimization problem.
In particular, a classical Newton method applied to the corresponding
first-order optimality conditions requires additional differentiability,
which may fail at degenerate parametrizations.
Even when the required derivatives exist along the iterates,
the resulting linear systems may become ill-conditioned as a trial
parametrization approaches elementwise degeneration.
For higher-order finite elements, this difficulty can arise even when
the degeneration is confined to an isolated point in the interior of
an element. These issues can complicate the numerical solution of
\eqref{eq:MMS-direct}.

To overcome this challenge, we replace $\cE_{\gamma}$ in \eqref{eq:MMS-direct} by an iteration-dependent \emph{modified energy} $\cE_{\mathrm{m},\gamma}^k$.

\begin{definition}
\label{def:modified-energy}
A functional
$\cE_{\mathrm{m},\gamma}^k : [V_h]^2 \to \R$
is called an \emph{admissible modified energy} for
$\cE_{\gamma}$ at $X_h^k$ if it satisfies the following conditions:
\begin{enumerate}
\item[\textnormal{(M1)}]
$\cE_{\mathrm{m},\gamma}^k(X_h^k) = \tfrac{1}{2}\cE_{\gamma}(X_h^k)$;
\item[\textnormal{(M2)}]
$\cE_{\mathrm{m},\gamma}^k(X) \geq \cE_{\gamma}(X) - \tfrac{1}{2}\cE_{\gamma}(X_h^k)$
for all $X \in [V_h]^2$;
\item[\textnormal{(M3)}]
The functional $\cE_{\mathrm{m},\gamma}^k$ is convex and
continuously differentiable on $[V_h]^2$.
\end{enumerate}
\end{definition}

Conditions (M1)--(M2) imply that
$$X \mapsto \cE_{\mathrm{m},\gamma}^k(X) + \tfrac12\cE_{\gamma}(X_h^k)$$ is majorizer of $\cE_{\gamma}$ that coincides with $\cE_{\gamma}(X)$ at
$X = X_h^k$. Thus, the modified energy is a convex surrogate anchored at the previous time level. These conditions also imply first-variation consistency whenever
$\cE_\gamma$ is differentiable at $X_h^k$. Indeed, the difference
\[
\cE_{\mathrm{m},\gamma}^k(X)
+
\frac12\cE_\gamma(X_h^k)
-
\cE_\gamma(X)
\]
is nonnegative and vanishes at $X_h^k$. Its first variation must therefore
vanish at that point. Hence,
\begin{equation}\label{eq:modified-energy-consistency}
D\cE_{\mathrm{m},\gamma}^k(X_h^k)
=
D\cE_\gamma(X_h^k).
\end{equation}

Given such a modified energy, we replace \eqref{eq:MMS-direct}
by
\begin{equation}\label{eq:MMS-framework}
X_h^{k+1} \in \mathop{\arg\min}_{X \in [V_h]^2}
\biggl\{ \frac{1}{2\tau}\|X - X_h^k\|_{L^2(\Gamma_h^k)}^2
+ \cE_{\mathrm{m},\gamma}^k(X)\biggr\}.
\end{equation}

\begin{theorem}
\label{thm:abstract}
Suppose that $\cE_{\mathrm{m},\gamma}^k$ is an admissible modified energy for
$\cE_{\gamma}$ at $X_h^k$. Assume that $|\p_\xi X_h^k| > 0$
on each element. Then the following statements hold.
\begin{enumerate}
\item[\textnormal{(i)}] Problem \eqref{eq:MMS-framework} has
a unique solution $X_h^{k+1} \in [V_h]^2$. The solution is characterized by
\begin{equation}\label{eq:abstract-EL}
\frac{1}{\tau}\int_\I (X_h^{k+1} - X_h^k)\cdot v_h\,
|\p_\xi X_h^k|\,\d\xi
+ D \cE_{\mathrm{m},\gamma}^k(X_h^{k+1})(v_h)= 0,
\,\,\,\forall\, v_h \in [V_h]^2.
\end{equation}

\item[\textnormal{(ii)}] The solution satisfies the one-step energy estimate
\begin{equation}\label{eq:abstract-dissipation}
\cE_{\gamma}(X_h^{k+1})
+ \frac{1}{2\tau}\|X_h^{k+1} - X_h^k\|_{L^2(\Gamma_h^k)}^2
\leq \cE_{\gamma}(X_h^k).
\end{equation}
\end{enumerate}
\end{theorem}

\Proof
(i)~The velocity penalty
$V(X) := \tfrac{1}{2\tau}\|X - X_h^k\|_{L^2(\Gamma_h^k)}^2$ is
strictly convex, and continuously differentiable on
$[V_h]^2$, since $|\p_\xi X_h^k| > 0$ guarantees that $\|\cdot\|_{L^2(\Gamma_h^k)}$ is a norm on $[V_h]^2$.
 Condition (M3) shows that  $\cE_{\mathrm{m},\gamma}^k$ is convex and continuously differentiable, too.
Their sum is therefore strictly convex,  coercive, and continuously differentiable on the finite-dimensional space $[V_h]^2$. Therefore, it has a unique minimizer, and its first-order optimality condition is \eqref{eq:abstract-EL}.

(ii)~To prove the energy estimate, we use $X_h^k$ as a competitor in \eqref{eq:MMS-framework}. The minimality of $X_h^{k+1}$ and condition
(M1) give
\[
\frac{1}{2\tau}\|X_h^{k+1} - X_h^k\|_{L^2(\Gamma_h^k)}^2
+ \cE_{\mathrm{m},\gamma}^k(X_h^{k+1})
\leq \cE_{\mathrm{m},\gamma}^k(X_h^k)
= \tfrac{1}{2}\cE_{\gamma}(X_h^k).
\]
Condition (M2), evaluated at $X_h^{k+1}$, gives
\[
\cE_{\mathrm{m},\gamma}^k(X_h^{k+1})
\geq
\cE_\gamma(X_h^{k+1})
-\frac12\cE_\gamma(X_h^k).
\]
Combining these two inequalities yields
\eqref{eq:abstract-dissipation}.
\qed

By \Cref{thm:abstract}, we can obtain the following statement: given an elementwise nondegenerate discrete parametrization $X_h^k$, the next iteration exists uniquely and satisfies \eqref{eq:abstract-dissipation} for any
$\tau>0$.  \Cref{thm:abstract} presents an energy-stable numerical scheme with arbitrary $\tau>0$ for simulating the $L^2$-gradient flow
of $\cE$ by the construction of an admissible modified energy satisfying (M1)--(M3). The theorem is conditional on the current parametrization: it does not prove that $X_h^{k+1}$ remains elementwise nondegenerate. Thus, the scheme can be advanced recursively only while the computed parametrizations remain elementwise nondegenerate.

The remainder of this section carries
out this construction in three concrete settings.

\subsection{Isotropic case: recovery of Dziuk's scheme}
\label{sec:mms-iso}
We begin with isotropic CSF:
$\p_t X = -\kappa N$, where the corresponding energy is the length functional  $$\mathcal{E}(X):=L(X)=\int_{\I} |\p_{\xi} X| \, \d \xi.$$
The construction relies on the elementary inequality
\begin{equation}\label{eq:alg-ineq}
\frac{a^2}{2b}
\geq
a-\frac{b}{2},
\qquad
a\geq0,\quad b>0,
\end{equation}
with equality if and only if $a=b$.

\begin{definition}
\label{def:iso-modified}
For $X, X_h^k \in [V_h]^2$ with $|\p_\xi X_h^k| > 0$ on each element, define
\begin{equation}\label{eq:Emix-iso}
\cE^k_{\mathrm{m}}(X)
:= \int_\I \frac{|\p_\xi X|^2}{2\,|\p_\xi X_h^k|}\,\d\xi.
\end{equation}
\end{definition}
Clearly, the first variation of $\cE^k_{\mathrm{m}}$ is
\begin{equation}\label{eq:iso-modified-first-variation}
    D\cE^k_{\mathrm{m}}(X)(v) =\l.\frac{\d }{\d \varepsilon}\r|_{ \varepsilon=0}\int_\I \frac{|\p_\xi X+\varepsilon \p_\xi v|^2}{2\,|\p_\xi X_h^k|}\,\d\xi=\int_\I \frac{\p_\xi X \cdot \p_\xi v}
{|\p_\xi X_h^k|}\,\d\xi.
\end{equation}

\begin{lemma}
\label{lem:iso-modified}
The functional $\cE_{\mathrm{m}}^k(X)$ defined by
\eqref{eq:Emix-iso} is an admissible modified energy for $\mathcal{E}$ at $X_h^k$.
\end{lemma}

\Proof Substitution of $X=X_h^k$ gives
\[
\cE_{\mathrm{m}}^k(X_h^k)
=
\frac12
\int_\I|\p_\xi X_h^k|\,\d\xi
=
\frac12\cE(X_h^k).
\]
Thus, condition (M1) holds.

We next apply \eqref{eq:alg-ineq} pointwise with
$a = |\p_\xi X|$, $b = |\p_\xi X_h^k|$. After integration over $\I$, we obtain
\[
\cE_{\mathrm{m}}^k(X)
\geq
\cE(X)-\frac12\cE(X_h^k),
\]
which is condition (M2).

Finally, the second variation is
\[
    D^2\cE_{\mathrm{m}}^k(X)(v,w)=\int_\I \frac{\p_\xi w \cdot \p_\xi v}
{|\p_\xi X_h^k|}\,\d\xi.
\]
In particular, it satisfies
$D^2\cE_{\mathrm{m}}^k(X)(v,v)\geq0$. Therefore, $\cE_{\mathrm{m}}^k$ is convex. Its quadratic form also shows that it is continuously differentiable on all of $[V_h]^2$. Hence,
condition (M3) holds.
\qed

Combining \Cref{lem:iso-modified} with \Cref{thm:abstract} gives the main result of this subsection.

\begin{theorem}[MMS for CSF]
\label{thm:mms-iso}
Assume $|\p_\xi X_h^k| > 0$ on each element.  Then the
minimization problem
\begin{equation}\label{eq:MMS-iso}
X_h^{k+1} \in \mathop{\arg\min}_{X \in [V_h]^2}
\biggl\{ \frac{1}{2\tau}\|X - X_h^k\|_{L^2(\Gamma_h^k)}^2
+ \cE_{\mathrm{m}}^k(X) \biggr\}
\end{equation}
has a unique solution, characterized by the Euler--Lagrange equation
\begin{equation}\label{eq:EL-Dziuk}
\frac{1}{\tau}\int_\I (X_h^{k+1} - X_h^k)\cdot v_h\,
|\p_\xi X_h^k|\,\d\xi
+ \int_\I \frac{\p_\xi X_h^{k+1} \cdot \p_\xi v_h}
{|\p_\xi X_h^k|}\,\d\xi = 0,
\qquad \forall\, v_h \in [V_h]^2.
\end{equation}
Moreover, for every $\tau>0$, the following one-step
energy estimate holds:
\begin{equation}\label{eq:energy-iso}
\cE(X_h^{k+1})
+ \frac{1}{2\tau}\|X_h^{k+1} - X_h^k\|_{L^2(\Gamma_h^k)}^2
\leq \cE(X_h^k).
\end{equation}
\end{theorem}

\begin{remark}\label{cor:dziuk}

For $n=1$, \eqref{eq:EL-Dziuk} recovers the consistent-mass form of
Dziuk's linearly implicit polygonal PFEM for CSF \cite{Dziuk1994};
for $n>1$, it defines the corresponding conforming piecewise-polynomial
Galerkin extension. In both cases, Theorem~\ref{thm:mms-iso}
characterizes the scheme as the Euler--Lagrange equation of a minimizing-movement step and provides a direct variational proof of the one-step
length estimate.
\end{remark}

\subsection{Weakly anisotropic case}
\label{sec:mms-weak}

We now extend the framework to weakly anisotropic CSF. The key technical ingredient is the following convexity result.

\begin{lemma}
\label{lem:gamma-sq-convex}
Let $\gamma\in C^2(\mathbb S^1)$ be positive and weakly anisotropic in the sense of \eqref{eq:weak-aniso}.
Then, $\widehat\gamma$ is convex
on $\R^2$. Moreover, the function
\[
\phi_\gamma(p):=\frac12\widehat\gamma(p)^2,
\qquad p\in\R^2,
\]
is convex and belongs to $C^1(\R^2)$, with
\begin{equation}\label{eq:A-gamma-extension}
\nabla\phi_\gamma(p):=\mathcal A_\gamma(p)
=
\begin{cases}
\widehat\gamma(p)\nabla\widehat\gamma(p),&p\neq0,\\
0,&p=0.
\end{cases}
\end{equation}
\end{lemma}

\Proof
Set
\[
M_0:=\max_{N\in\mathbb S^1}\gamma(N).
\]
Then \(0\leq\widehat\gamma(p)\leq M_0|p|\), so
\(\widehat\gamma\) is continuous at the origin. For \(p\neq0\), set $N_p:=\frac{p}{|p|}$. Positive one-homogeneity and Euler's identity yield
\[
\nabla^2\widehat\gamma(p)p=0,
\qquad
\nabla^2\widehat\gamma(p)
=\frac{1}{|p|}\nabla^2\widehat\gamma(N_p).
\]
For any $z\in \R^2$, define
\[z_\perp:=z-(z\cdot N_p)N_p.\]
Since $z_\perp\perp N_p$, the symmetry of the Hessian and
\eqref{eq:weak-aniso} imply
\[
\nabla^2\widehat\gamma(p)z\cdot z
=
\frac{1}{|p|}
\nabla^2\widehat\gamma(N_p)z_\perp\cdot z_\perp
\geq0.
\]
Thus
\(\nabla^2\widehat\gamma\) is positive semidefinite on $\R^2\setminus\{0\}$,
and the restriction of $\widehat\gamma$ to any affine line avoiding
the origin is convex.

For a line through the
origin, positive homogeneity gives, for \(d\neq0\),
\[
\widehat\gamma(td)
=
\begin{cases}
t\,\widehat\gamma(d),&t\geq0,\\
(-t)\,\widehat\gamma(-d),&t<0.
\end{cases}
\]
This piecewise-linear function is continuous at $t=0$, and its
left and right slopes satisfy
\[
-\widehat\gamma(-d)<\widehat\gamma(d)
\]
by the positivity of $\widehat\gamma$.
Hence this restriction is also convex.
Since its restriction to every affine line is convex,
$\widehat\gamma$ is convex on $\R^2$.

Since \(\widehat\gamma\geq0\) and \(r\mapsto r^2/2\) is convex and
nondecreasing on \([0,\infty)\), \(\phi_\gamma\) is convex. Moreover,
\[
0\leq\phi_\gamma(p)
\leq\frac{M_0^2}{2}|p|^2,
\]
so \(\phi_\gamma\) is Fr\'echet differentiable at the origin with
\(\nabla\phi_\gamma(0)=0\). For \(p\neq0\), the chain rule gives
\[
\nabla\phi_\gamma(p)
=\widehat\gamma(p)\nabla\widehat\gamma(p).
\]
Finally, since \(\nabla\widehat\gamma\) is zero-homogeneous, with
\[
M_1:=\max_{N\in\mathbb S^1}
|\nabla\widehat\gamma(N)|<\infty,
\]
we obtain
\[
|\nabla\phi_\gamma(p)|
\leq M_0M_1|p|
\longrightarrow0
\qquad\text{as }p\to0.
\]
Thus \(\phi_\gamma\in C^1(\R^2)\), and
\eqref{eq:A-gamma-extension} follows.
\qed

\begin{definition}
\label{def:aniso-modified}
For $X, X_h^k\in [V_h]^2$ with $|\p_\xi X_h^k| > 0$ on
each element, we define
\begin{equation}\label{eq:Emix-aniso}
\cE_{\mathrm{m},\gamma}^k(X)
:= \int_\I
\frac{\widehat\gamma^2\bigl(p\bigr)}
     {2\,\widehat\gamma(p^{k})}\,\d\xi,\qquad p^{k}:= (\p_\xi X_h^{k})^\perp,\quad p:=p(X)=(\p_\xi X)^\perp.
\end{equation}
\end{definition}

Because $\widehat\gamma(p^k)>0$ on every element, the denominator in
\eqref{eq:Emix-aniso} is well defined. The globally defined map
$\mathcal A_\gamma$ in \eqref{eq:A-gamma-extension} gives
\begin{equation}\label{eq:aniso-modified-first-variation}
D\cE_{\mathrm{m},\gamma}^k(X)(v)
=
\int_\I
\frac{
\mathcal A_\gamma\bigl(p\bigr)
\cdot(\p_\xi v)^\perp
}{
\widehat\gamma(p^k)
}
\,\d\xi.
\end{equation}

\begin{lemma}
\label{lem:aniso-modified}
Assume that $\gamma$ is weakly anisotropic and that
$|\p_\xi X_h^k| > 0$ on each element. Then $\cE_{\mathrm{m},\gamma}^k$ is an
admissible modified energy for $\cE_\gamma$ at $X_h^k$.
\end{lemma}

\Proof
Substitution of $X=X_h^k$ into \eqref{eq:Emix-aniso} gives
\[
\cE_{\mathrm{m},\gamma}^k(X_h^k)
=
\frac12
\int_\I\widehat\gamma(p^k)\,\d\xi
=
\frac12\cE_\gamma(X_h^k).
\]
Thus, condition (M1) holds.

We apply \eqref{eq:alg-ineq} pointwise with
$a=\widehat\gamma\bigl(p\bigr)\geq0$, $b=\widehat\gamma(p^k)>0$.
After integration, we obtain
\[
\cE_{\mathrm{m},\gamma}^k(X)
\geq
\cE_\gamma(X)
-\frac12\cE_\gamma(X_h^k).
\]
This proves condition (M2).

By Lemma~\ref{lem:gamma-sq-convex}, the mapping
$p\longmapsto\frac12\widehat\gamma(p)^2$
is convex and continuously differentiable on $\R^2$. The denominator in \eqref{eq:Emix-aniso} is positive and fixed on each element. Since $X\mapsto p(X)$ is linear, $\cE_{\mathrm{m},\gamma}^k$ is a convex $C^1$ functional on
$[V_h]^2$. This conclusion also applies to trial functions with collapsed
edges. Hence, condition (M3) holds. \qed

\begin{theorem}[MMS for weakly anisotropic CSF]
\label{thm:mms-weak}
Assume that $\gamma$ is weakly anisotropic and that $|\p_\xi X_h^k| > 0$ on every element. Then the minimization problem
\begin{equation}\label{eq:MMS-aniso-weak}
X_h^{k+1} \in \mathop{\arg\min}_{X \in [V_h]^2}
\biggl\{ \frac{1}{2\tau}\|X - X_h^k\|_{L^2(\Gamma_h^k)}^2
+ \cE_{\mathrm{m},\gamma}^k(X) \biggr\}
\end{equation}
has a unique solution, characterized by the Euler--Lagrange equation
\begin{equation}\label{eq:EL-aniso-weak}
\frac{1}{\tau}\int_\I (X_h^{k+1} - X_h^k)\cdot v_h\,
|\p_\xi X_h^k|\,\d\xi
+ \int_\I
{\frac{\mathcal A_\gamma(p^{k+1})
\cdot (\p_\xi v_h)^\perp}
{\widehat\gamma(p^k)}}\,\d\xi = 0, \quad \forall v_h\in[V_h]^2.
\end{equation}
Moreover, for every $\tau>0$, the following one-step energy estimate holds:
\begin{equation}\label{eq:energy-aniso-weak}
\cE_\gamma(X_h^{k+1})
+ \frac{1}{2\tau}\|X_h^{k+1} - X_h^k\|_{L^2(\Gamma_h^k)}^2
\leq \cE_\gamma(X_h^k).
\end{equation}
\end{theorem}

\Proof
The result follows from Lemma~\ref{lem:aniso-modified} and
Theorem~\ref{thm:abstract}. Equation~\eqref{eq:EL-aniso-weak} follows
from \eqref{eq:aniso-modified-first-variation}.
\qed

\subsection{Strongly anisotropic case via convex splitting}
\label{sec:mms-strong}

When $\gamma$ is strongly anisotropic, the functional $\cE_{\mathrm{m},\gamma}^k(X)$ in \Cref{def:aniso-modified}
need not be convex. Therefore, condition (M3) may fail, and the preceding construction cannot be applied directly. We
restore the required convexity by combining the modified-energy framework
with convex splitting. Related splitting techniques have been used to construct energy-stable methods for gradient flows with nonconvex energies \cite{Eyre1998,Shen2010,dai2022convex}. Here, we adapt this idea to a parametric finite element discretization.

\begin{definition}
A decomposition $\gamma = \gamma_c - \gamma_e$ is called a
\emph{convex splitting} if both $\gamma_c$ and $\gamma_e$ are positive and weakly anisotropic. The anisotropic energy is then splitted correspondingly into $\cE_{\gamma}=\cE_{\gamma_c}- \cE_{\gamma_e}$.
\end{definition}

\begin{remark}
\label{rem:splitting}
Each positive $\gamma\in C^2(\S^1)$ admits such a splitting when $\gamma$ is strongly anisotropic. Set
\[C_{\min}
:=\max\left\{0, -\min_{N \in \mathbb{S}^1}
\nabla^2 \,\widehat\gamma(N)\,T\cdot T\right\}, \]
where $T$ is a unit tangent vector orthogonal to $N$.
For a strongly anisotropic density, $C_{\min}>0$. Any choice $C\geq C_{\min}$ gives \[
\gamma_c=\gamma+C,
\qquad
\gamma_e=C.
\]
Indeed,
\[
\nabla^2\,\widehat\gamma_c(N)\,T\cdot T
= \nabla^2\,\widehat\gamma(N)\,T\cdot T + C \geq 0,
\qquad
\nabla^2\,\widehat\gamma_e(N)\,T\cdot T
=C>0.
\]
Thus, both components are weakly anisotropic. For the $q$-fold density
$\gamma(N)=1+\beta\cos(q\theta)$,
$C_{\min}=\max\left\{
0,\, \beta(q^2-1)-1
\right\}$. The nonnegative quantity $C-C_{\min}$ measures the amount of
over-splitting beyond the minimum required to make $\gamma_c$ weakly
anisotropic. Its numerical effect is examined in \Cref{sec:strong-aniso-conv}.
\end{remark}

\begin{definition}
\label{def:cs-modified}
Let $\gamma=\gamma_c-\gamma_e$ be a convex splitting. For
$X, X_h^k\in[V_h]^2$ with $|\p_\xi X_h^k| > 0$, define
\begin{equation}\label{eq:cs-modified}
\cE_{\mathrm{mc},\gamma}^k(X)
:=\cE_{\mathrm{m},\gamma_c}^k(X)
-D\cE_{\gamma_e}(X_h^k)(X-X_h^k)-
\frac12\cE_{\gamma_e}(X_h^k),
\end{equation}
where $\cE_{\mathrm{m},\gamma_c}^k$ is defined by
\eqref{eq:Emix-aniso} with $\gamma$ replaced by $\gamma_c$. The
linearized explicit contribution is
\begin{equation}\label{eq:first-var}
D\cE_{\gamma_e}(X_h^k)(X-X_h^k)
=\int_\I
\nabla\widehat\gamma_e(p^k)
\cdot
\bigl(p(X)-p^k\bigr)
\,\d\xi.
\end{equation}
\end{definition}

\begin{remark}
The constant $-\tfrac12\cE_{\gamma_e}(X_h^k)$ in
\eqref{eq:cs-modified} does not affect the minimizer. It is included
to ensure that condition (M1) holds for the complete energy $\cE_\gamma$.
\end{remark}

\begin{lemma}
\label{lem:cs-modified}
Suppose that $\gamma = \gamma_c - \gamma_e$ is a convex splitting and that
$|\p_\xi X_h^k| > 0$ on each element. Then
$\cE_{\mathrm{mc},\gamma}^{k}(X)$ is an admissible modified energy for $\cE_\gamma$ at $X_h^k$.
\end{lemma}

\Proof At $X=X_h^k$, the linear term in \eqref{eq:cs-modified} vanishes.
Lemma~\ref{lem:aniso-modified} applied to $\gamma_c$, gives
\[
\cE_{\mathrm{m},\gamma_c}^k(X_h^k)
=
\frac12\cE_{\gamma_c}(X_h^k).
\]
Therefore,
\[
\cE_{\mathrm{mc},\gamma}^k(X_h^k)
=
\frac12\cE_{\gamma_c}(X_h^k)
-
\frac12\cE_{\gamma_e}(X_h^k)
=
\frac12\cE_\gamma(X_h^k).
\]
This proves condition (M1).

Because $\gamma_e$ is weakly anisotropic,
Lemma~\ref{lem:gamma-sq-convex} shows that
$\widehat\gamma_e$ is convex. The supporting-hyperplane inequality
therefore gives
\[
\cE_{\gamma_e}(X)
\geq
\cE_{\gamma_e}(X_h^k)
+
D\cE_{\gamma_e}(X_h^k)(X-X_h^k).
\]
Equivalently,
\[
-
D\cE_{\gamma_e}(X_h^k)(X-X_h^k)
\geq
-\cE_{\gamma_e}(X)
+
\cE_{\gamma_e}(X_h^k).
\]
Lemma~\ref{lem:aniso-modified}, applied to $\gamma_c$, also gives
\[
\cE_{\mathrm{m},\gamma_c}^k(X)
\geq
\cE_{\gamma_c}(X)
-
\frac12\cE_{\gamma_c}(X_h^k).
\]
Adding these inequalities and subtracting
$\frac12\cE_{\gamma_e}(X_h^k)$ yields
\[
\cE_{\mathrm{mc},\gamma}^k(X)
\geq
\cE_\gamma(X)
-
\frac12\cE_\gamma(X_h^k).
\]
Thus, condition (M2) holds.

Finally, $\cE_{\mathrm{m},\gamma_c}^k$ is convex and $C^1$. The
remaining $X$-dependent term in \eqref{eq:cs-modified} is affine.
Therefore, $\cE_{\mathrm{mc},\gamma}^k$ is convex and $C^1$ on
$[V_h]^2$. This proves condition (M3).
\qed

\begin{theorem}[MMS for strongly anisotropic CSF] \label{thm:mms-strong}
Suppose that $\gamma=\gamma_c-\gamma_e$ is a convex splitting and that
$|\p_\xi X_h^k|>0$ on every element. Then
\begin{equation}\label{eq:MMS-aniso-strong}
X_h^{k+1}
\in
\mathop{\arg\min}_{X\in[V_h]^2}
\left\{
\frac{1}{2\tau}
\|X-X_h^k\|_{L^2(\Gamma_h^k)}^2
+
\cE_{\mathrm{mc},\gamma}^k(X)
\right\}
\end{equation}
has a unique solution. This solution satisfies
\begin{equation}\label{eq:EL-aniso-strong}
\begin{aligned}
&
\frac{1}{\tau}
\int_\I
(X_h^{k+1}-X_h^k)\cdot v_h
|\p_\xi X_h^k|\,\d\xi
+
\int_\I
\frac{
\mathcal A_{\gamma_c}(p^{k+1})
\cdot(\p_\xi v_h)^\perp
}{
\widehat\gamma_c(p^k)
}
\,\d\xi
\\
&\quad
-
\int_\I
\nabla\widehat\gamma_e(p^k)
\cdot(\p_\xi v_h)^\perp
\,\d\xi
=0,
\qquad
\forall v_h\in[V_h]^2.
\end{aligned}
\end{equation}
Moreover, for every $\tau>0$,
\begin{equation}\label{eq:energy-aniso-strong}
\cE_\gamma(X_h^{k+1})
+
\frac{1}{2\tau}
\|X_h^{k+1}-X_h^k\|_{L^2(\Gamma_h^k)}^2
\leq
\cE_\gamma(X_h^k).
\end{equation}
\end{theorem}

\Proof Lemma~\ref{lem:cs-modified} and Theorem~\ref{thm:abstract} give unique
solvability and the energy estimate. Differentiating
\eqref{eq:cs-modified} with respect to $X$ gives
\eqref{eq:EL-aniso-strong}. The constant term has zero derivative.
\qed

\section{Extension to area-preserving curve shortening flow}
\label{sec:ap-csf}
We now incorporate an area constraint into the minimizing-movement
framework through a single formulation for isotropic, weakly
anisotropic, and strongly anisotropic area-preserving curve shortening
flows (AP-CSFs). The construction retains the appropriate modified
energy from the preceding section and adds only one scalar constraint.
Consequently, every global minimizer preserves the discrete signed area
and satisfies an energy dissipation inequality.

The area-preserving anisotropic CSF is the constrained
$L^2$-gradient flow of $\cE_\gamma$; the classical isotropic AP-CSF is
obtained by taking $\gamma\equiv1$. For $X\in[V_h]^2$, we define the
algebraic signed area enclosed by the discrete parametric curve
$\Gamma(X)$ as
\[
A(X)
:= \frac{1}{2}\int_\I
X\cdot(\p_\xi X)^\perp\,\d\xi.
\]
This functional is quadratic but indefinite with respect to the finite element coefficient vector.
Let $A_0 := A(X_h^0)\neq0$ be the initial signed area, and define the constraint set
\[
\cC_{A_0}:=\bigl\{\, X \in [V_h]^2 :
A(X) = A_0\,\bigr\}.
\]

Throughout this section, we assume that
$X_h^k\in\cC_{A_0}$ and that
$|\p_\xi X_h^k|>0$ on every element. Periodic integration by parts
gives
\begin{equation}\label{DA}
D A(X)(v)
=\frac{1}{2}\int_\I
\left[
v\cdot(\p_\xi X)^\perp
+X\cdot(\p_\xi v)^\perp
\right]\d\xi =
\int_\I
(\p_\xi X)^\perp\cdot v\,\d\xi.
\end{equation}

\begin{definition}
\label{def:MMS-AP}
Let $\gamma\in C^2(\mathbb S^1)$ be positive. We define the unified modified energy by
\[
\widetilde\cE_{\mathrm{m},\gamma}^k(X)
:=
\begin{cases}
\cE_{\mathrm{m},\gamma}^k(X),
& \gamma\ \text{is weakly anisotropic},\\[1mm]
\cE_{\mathrm{mc},\gamma}^k(X),
& \gamma\ \text{is strongly anisotropic}.
\end{cases}
\]
The two modified energies are defined in
\Cref{def:aniso-modified,def:cs-modified}, respectively. In the
strongly anisotropic case, one may use the constant convex splitting
from \Cref{rem:splitting} with $C\geq C_{\min}$.
\end{definition}

Given $X_h^k\in\cC_{A_0}$, the area-preserving minimizing-movement step seeks
\begin{equation}\label{eq:MMS-AP-aniso}
X_h^{k+1}
\in \mathop{\arg\min}_{X\in\cC_{A_0}}
\biggl\{
\frac{1}{2\tau}\|X-X_h^k\|_{L^2(\Gamma_h^k)}^2
+\widetilde\cE_{\mathrm{m},\gamma}^k(X)
\biggr\}.
\end{equation}

\begin{theorem}[MMS for the AP-CSF]
\label{thm:ap-aniso}
Let $\gamma\in C^2(\mathbb S^1)$ be positive, and choose
$\widetilde\cE_{\mathrm{m},\gamma}^k$ according to
\Cref{def:MMS-AP}. Assume that $\tau>0$,
$X_h^k\in\cC_{A_0}$, and $|\p_\xi X_h^k|>0$ on every element.
Then the following statements hold.
\begin{enumerate}
\item[\textnormal{(i)}]
Problem \eqref{eq:MMS-AP-aniso} admits at least one global minimizer
$X_h^{k+1}\in\cC_{A_0}$.

\item[\textnormal{(ii)}]
Every global minimizer preserves the signed area,
$A(X_h^{k+1})=A_0$, and satisfies the discrete energy dissipation
inequality
\begin{equation}\label{eq:energy-AP-aniso}
\cE_\gamma(X_h^{k+1})
+\frac{1}{2\tau}
\|X_h^{k+1}-X_h^k\|_{L^2(\Gamma_h^k)}^2
\leq \cE_\gamma(X_h^k).
\end{equation}

\item[\textnormal{(iii)}]
Every global minimizer satisfies the constraint qualification
$D A(X_h^{k+1})\neq0$. Moreover, there exists a Lagrange multiplier
$\lambda^{k+1}\in\R$ associated with this minimizer such that
\begin{equation}
\begin{aligned}
&
\frac{1}{\tau}
\int_\I
(X_h^{k+1}-X_h^k)\cdot v_h
|\p_\xi X_h^k|\,\d\xi
+
D\widetilde\cE_{\mathrm{m},\gamma}^k
(X_h^{k+1})(v_h)
\\
&\qquad
+
\lambda^{k+1}
\int_\I
v_h\cdot(\p_\xi X_h^{k+1})^\perp\,\d\xi
=0,
\qquad
\forall\,v_h\in[V_h]^2,
\end{aligned}
\label{eq:EL-AP-aniso}
\end{equation}
together with $A(X_h^{k+1})=A_0$. The derivative of the modified
energy in \eqref{eq:EL-AP-aniso} is given by the energy-variation
terms in \eqref{eq:EL-aniso-weak} or
\eqref{eq:EL-aniso-strong}, according to the type of anisotropy.
\end{enumerate}
\end{theorem}

\Proof By \Cref{lem:aniso-modified,lem:cs-modified}, the selected functional
$\widetilde\cE_{\mathrm{m},\gamma}^k$ is an admissible modified energy
for $\cE_\gamma$ at $X_h^k$. Define
\[
F_h^k(X)
:=
\frac{1}{2\tau}
\|X-X_h^k\|_{L^2(\Gamma_h^k)}^2
+
\widetilde\cE_{\mathrm{m},\gamma}^k(X).
\]
Property (M3) implies that $F_h^k$ is continuous. Moreover, (M2) and the positivity of $\cE_\gamma$ give
\[
F_h^k(X)
\geq
\frac{1}{2\tau}
\|X-X_h^k\|_{L^2(\Gamma_h^k)}^2
-\frac12\cE_\gamma(X_h^k).
\]
Since $|\p_\xi X_h^k|>0$ on every element,
$\|\cdot\|_{L^2(\Gamma_h^k)}$ is a norm on $[V_h]^2$. Hence, $F_h^k$ is coercive.
The set $\cC_{A_0}$ is closed and nonempty because
$A$ is continuous and $X_h^k\in\cC_{A_0}$. The finite-dimensional direct method therefore
gives a global minimizer of \eqref{eq:MMS-AP-aniso}.

Since $X_h^k$ is feasible, minimality and (M1) yield
\[
\frac{1}{2\tau}
\|X_h^{k+1}-X_h^k\|_{L^2(\Gamma_h^k)}^2
+
\widetilde\cE_{\mathrm{m},\gamma}^k(X_h^{k+1})
\leq
\frac12\cE_\gamma(X_h^k).
\]
Combining this inequality with (M2), evaluated at $X_h^{k+1}$,
proves \eqref{eq:energy-AP-aniso}. Membership in $\cC_{A_0}$ gives
$A(X_h^{k+1})=A_0$.

Finally, \eqref{DA} and $A_0\neq0$ imply
\[
D A(X_h^{k+1})(X_h^{k+1})
=2A(X_h^{k+1})
=2A_0
\neq0.
\]
Thus, the constraint qualification holds. Since the objective and
the constraint are continuously differentiable, the
finite-dimensional Lagrange multiplier theorem yields
\eqref{eq:EL-AP-aniso} together with the area constraint.
\qed

\begin{remark}
The isotropic case is included by taking
$\gamma\equiv1$. In this
case,
$\widehat\gamma(p)=|p|$, $\mathcal A_\gamma(p)=p$, and
\[
\widetilde\cE_{\mathrm{m},\gamma}^k(X)
=\int_\I
\frac{|\p_\xi X|^2}{2|\p_\xi X_h^k|}\,\d\xi
=\cE_{\mathrm{m}}^k(X).
\]
Thus, \eqref{eq:MMS-AP-aniso} reduces to the area-preserving isotropic minimizing-movement scheme. Moreover,
\[
D\widetilde\cE_{\mathrm{m},\gamma}^k
(X_h^{k+1})(v_h)
=
\int_\I
\frac{\p_\xi X_h^{k+1}\cdot\p_\xi v_h}
{|\p_\xi X_h^k|}\,\d\xi.
\]
Hence, \eqref{eq:EL-AP-aniso} is obtained from
\eqref{eq:EL-Dziuk} by adding one scalar area-multiplier term. The energy estimate \eqref{eq:energy-AP-aniso} then becomes the
corresponding length dissipation inequality.
\end{remark}

\begin{remark}
With the sign convention used in the preceding sections, the smooth
isotropic AP-CSF prescribes the normal velocity
\[
V_n=\overline{\kappa}-\kappa,
\qquad
\overline{\kappa}
:=\frac{1}{L(\Gamma)}
\int_\Gamma\kappa\,\d s.
\]
Because \eqref{eq:EL-AP-aniso} contains
$+\lambda^{k+1}D A(X_h^{k+1})$, the quantity
$-\lambda^{k+1}$ is the discrete counterpart of the nonlocal term
$\overline\kappa$.
An exact global minimizer preserves the
signed area exactly, whereas a numerically computed  solution satisfies the constraint only up to the nonlinear-solver tolerance. We refer the
reader to \cite{Jiang2023,Jiang2025,pei2023structure} for related
structure-preserving PFEMs for area-preserving curvature flows.
Lagrange-multiplier formulations for geometric flows are discussed in \cite{garcke2025structure}.
\end{remark}

\section{Numerical results}\label{sec:numerics}
This section presents numerical experiments for the proposed
MMS--PFEM schemes. The experiments examine three aspects of the
methods. First, we study the coupled space--time convergence behavior in the isotropic and weakly anisotropic cases. Second, we investigate
the mesh-dependent oscillations that arise for strong anisotropies.
Finally, we examine the discrete energy decay, area preservation, mesh
quality, and nonlinear-solver performance. The nonlinear
systems associated with
\eqref{eq:EL-Dziuk},
\eqref{eq:EL-aniso-weak},
\eqref{eq:EL-aniso-strong},
and the unified constrained system
\eqref{eq:EL-AP-aniso}
are solved by Newton's method, which is
initialized with the previous solution, $
X_h^{k+1,0}:=X_h^k$ at each time step. The iteration stops when the $\ell^\infty$-norm of the assembled
residual is below $10^{-12}$. All computations are performed in Firedrake using continuous piecewise-linear finite elements ($n=1$), and the codes are deposited in Github at https://github.com/Ganghui-Zhang/-L-2--MMS.
%For the AP-CSF schemes, this residual includes the scalar area-constraint equation.
\subsection{Convergence tests for isotropic and weakly anisotropic flows}

We first study the observed convergence behavior of the CSF and
AP-CSF schemes. The initial curve is the unit circle for the CSF and
the ellipse $x^2/4+y^2=1$ for the AP-CSF. Let $J$ denote the number of mesh elements and set
$h=J^{-1}$. We compute solutions using
$J=16,32,64,128,256,512$, and choose the time step as $\tau=h^2$.

Exact solutions are unavailable for most of these experiments.
Therefore, we compare solutions on successive mesh levels using the
manifold distance. For two simple closed curves $\Gamma_1$ and
$\Gamma_2$, this distance is defined by \cite{Zhao2021energy,Jiang2024second}
\begin{equation}
\mathrm{M}(\Gamma_1,\Gamma_2)
:=
|\Omega_1\mathbin{\triangle}\Omega_2|
=
|\Omega_1|+|\Omega_2|-2|\Omega_1\cap\Omega_2|,
\label{eq:manifold-distance}
\end{equation}
where $\Omega_i$ is the region enclosed by $\Gamma_i$ and
$\triangle$ denotes the symmetric difference.

\begin{figure}[h!]
    \centering
    \begin{tikzpicture}
        \pgfplotsset{
            every axis/.append style={
                width=0.47\linewidth, height=0.38\linewidth,
                xmode=log, ymode=log,
                xlabel style={font=\normalsize, yshift=2pt},
                ylabel style={font=\normalsize, yshift=-2pt},
                title style={font=\small},
                tick label style={font=\scriptsize},
                grid=both,
                grid style={dotted, gray!40},
                every axis plot/.append style={very thick},
                legend style={font=\small, draw=none, fill opacity=0.7,
                              inner sep=2pt, row sep=-2pt},
                legend cell align=left,
            },
        }
        \begin{groupplot}[
            group style={group size=2 by 2,
                         horizontal sep=1.8cm, vertical sep=2.2cm},
        ]
        % --- (a) Case 1: isotropic ---
        \nextgroupplot[
            title={(a) Case 1: isotropic, $\gamma \equiv 1$},
            xlabel={$h$}, ylabel={$E_{\mathrm{M}}$},
            legend pos=south east,
        ]
            \addplot[seabornblue!80!black, mark=*]
                table[x index=0, y index=1] {figure/eoc_iso_circle.dat};
                \addlegendentry{$\mathrm{T}=0.05$}
            \addplot[seabornred!80!black, mark=square*]
                table[x index=0, y index=2] {figure/eoc_iso_circle.dat};
                \addlegendentry{$\mathrm{T}=0.1$}
            \addplot[seaborngreen!80!black, mark=triangle*]
                table[x index=0, y index=3] {figure/eoc_iso_circle.dat};
                \addlegendentry{$\mathrm{T}=0.3$}
            \addplot[domain=0.003:0.08, dashed, black] {4*x^2};
                \addlegendentry{slope 2}
        % --- (b) Case 2: q-fold ---
        \nextgroupplot[
            title={(b) Case 2: $q$-fold, $(q,\beta)=(4,0.03)$},
            xlabel={$h$}, ylabel={$E_{\mathrm{M}}$},
            legend pos=south east,
        ]
            \addplot[seabornblue!80!black, mark=*]
                table[x index=0, y index=1] {figure/eoc_kfold_k4_weak_circle.dat};
        %        \addlegendentry{$\mathrm{T}=0.05$}
            \addplot[seabornred!80!black, mark=square*]
                table[x index=0, y index=2] {figure/eoc_kfold_k4_weak_circle.dat};
        %        \addlegendentry{$\mathrm{T}=0.1$}
            \addplot[seaborngreen!80!black, mark=triangle*]
                table[x index=0, y index=3] {figure/eoc_kfold_k4_weak_circle.dat};
        %        \addlegendentry{$\mathrm{T}=0.3$}
            \addplot[domain=0.003:0.08, dashed, black] {4*x^2};
       %         \addlegendentry{slope 2}
        % --- (c) Case 3: Riemannian ---
        \nextgroupplot[
            title={(c) Case 3: Riemannian, $M=1$, $G=\mathrm{diag}(1,2)$},
            xlabel={$h$}, ylabel={$E_{\mathrm{M}}$},
            legend pos=south east,
        ]
            \addplot[seabornblue!80!black, mark=*]
                table[x index=0, y index=1] {figure/eoc_riemannian_circle.dat};
         %       \addlegendentry{$\mathrm{T}=0.05$}
            \addplot[seabornred!80!black, mark=square*]
                table[x index=0, y index=2] {figure/eoc_riemannian_circle.dat};
        %        \addlegendentry{$\mathrm{T}=0.1$}
            \addplot[seaborngreen!80!black, mark=triangle*]
                table[x index=0, y index=3] {figure/eoc_riemannian_circle.dat};
       %         \addlegendentry{$\mathrm{T}=0.3$}
            \addplot[domain=0.003:0.08, dashed, black] {4*x^2};
      %          \addlegendentry{slope 2}
        % --- (d) Case 4: l^r-norm metric ---
        \nextgroupplot[
            title={(d) Case 4: $\ell^r$-norm metric, $r=4$},
            xlabel={$h$}, ylabel={$E_{\mathrm{M}}$},
            legend pos=south east,
        ]
            \addplot[seabornblue!80!black, mark=*]
                table[x index=0, y index=1] {figure/eoc_l4_circle.dat};
     %           \addlegendentry{$\mathrm{T}=0.05$}
            \addplot[seabornred!80!black, mark=square*]
                table[x index=0, y index=2] {figure/eoc_l4_circle.dat};
     %           \addlegendentry{$\mathrm{T}=0.1$}
            \addplot[seaborngreen!80!black, mark=triangle*]
                table[x index=0, y index=3] {figure/eoc_l4_circle.dat};
     %           \addlegendentry{$\mathrm{T}=0.3$}
            \addplot[domain=0.003:0.08, dashed, black] {4*x^2};
    %            \addlegendentry{slope 2}
        \end{groupplot}
    \end{tikzpicture}
    \caption{Adjacent-grid manifold distance
    $E_{\mathrm M}$ for the CSF. The initial curve is the unit circle. The results are shown at $\mathrm{T}=0.05$, $0.1$, and $0.3$
    for the isotropic density and three weakly anisotropic densities. The dashed lines have slope two.}
    \label{fig:eoc_csf}
\end{figure}
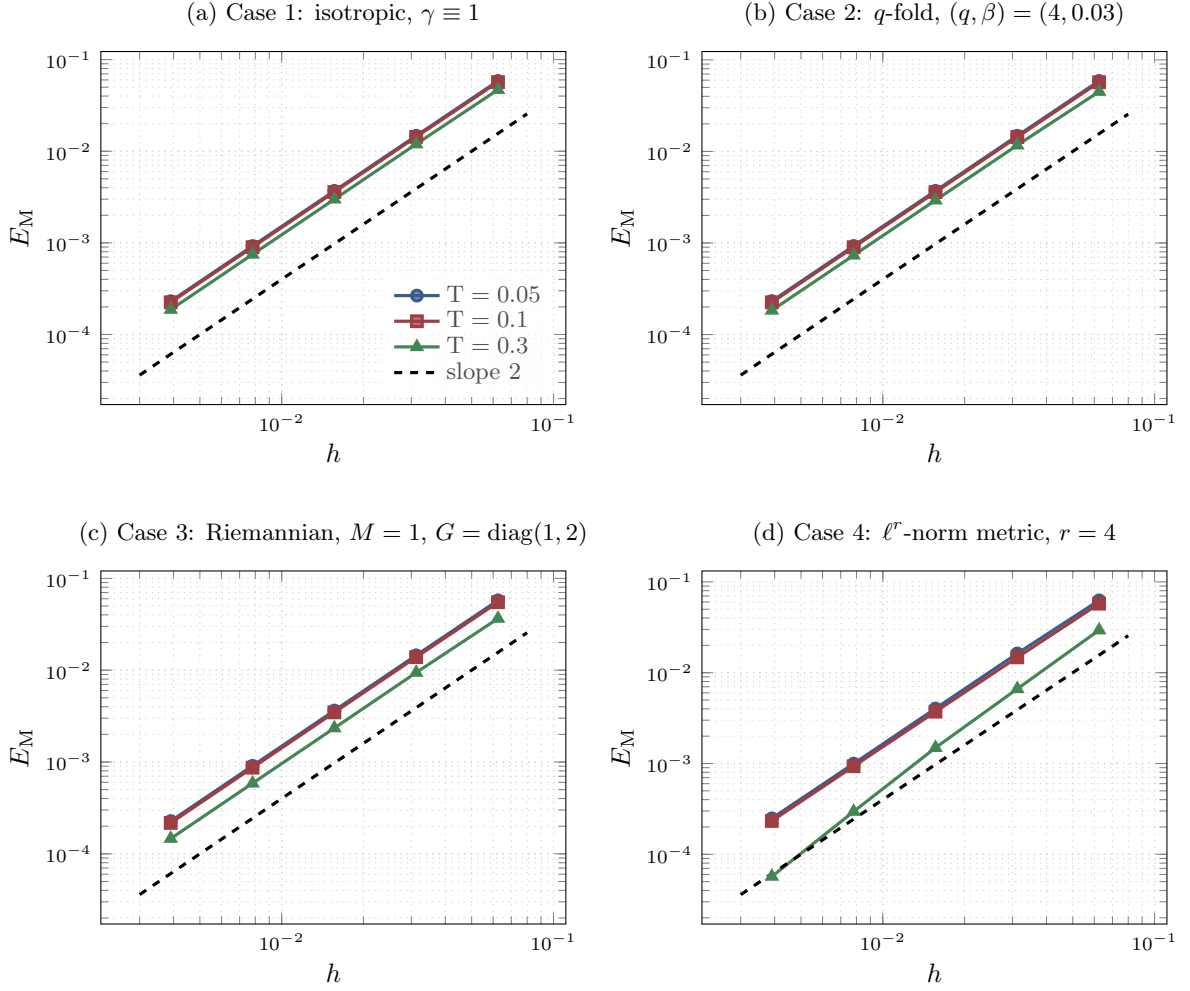

Let $\Gamma_{h,\tau}(\mathrm{T})$ denote the numerical curve at time $\mathrm{T}$.
The adjacent-grid error is
\begin{equation}
E_{\mathrm{M}}(h,\tau;\mathrm{T})
:=
\mathrm{M}\left(
\Gamma_{h,\tau}(\mathrm{T}),
\Gamma_{h/2,\tau/4}(\mathrm{T})
\right).
\label{eq:coupled-self-error}
\end{equation}
The backward Euler method is first order in time. Thus, the choice
$\tau=h^2$ balances an $O(\tau)$ temporal contribution with a
putative $O(h^2)$ spatial contribution. The resulting errors are
coupled space--time errors. Therefore, these experiments do not
separately identify the spatial and temporal convergence orders. We calculate the adjacent-grid order from
\begin{equation}
\operatorname{Order}(h,\tau;\mathrm{T})
:=
\frac{
\log\left(
E_{\mathrm{M}}(h,\tau;\mathrm{T})/
E_{\mathrm{M}}(h/2,\tau/4;\mathrm{T})
\right)}
{\log 2}.
\label{eq:observed-order}
\end{equation}
An order close to two is consistent with the expected combined error
when $\tau=h^2$.

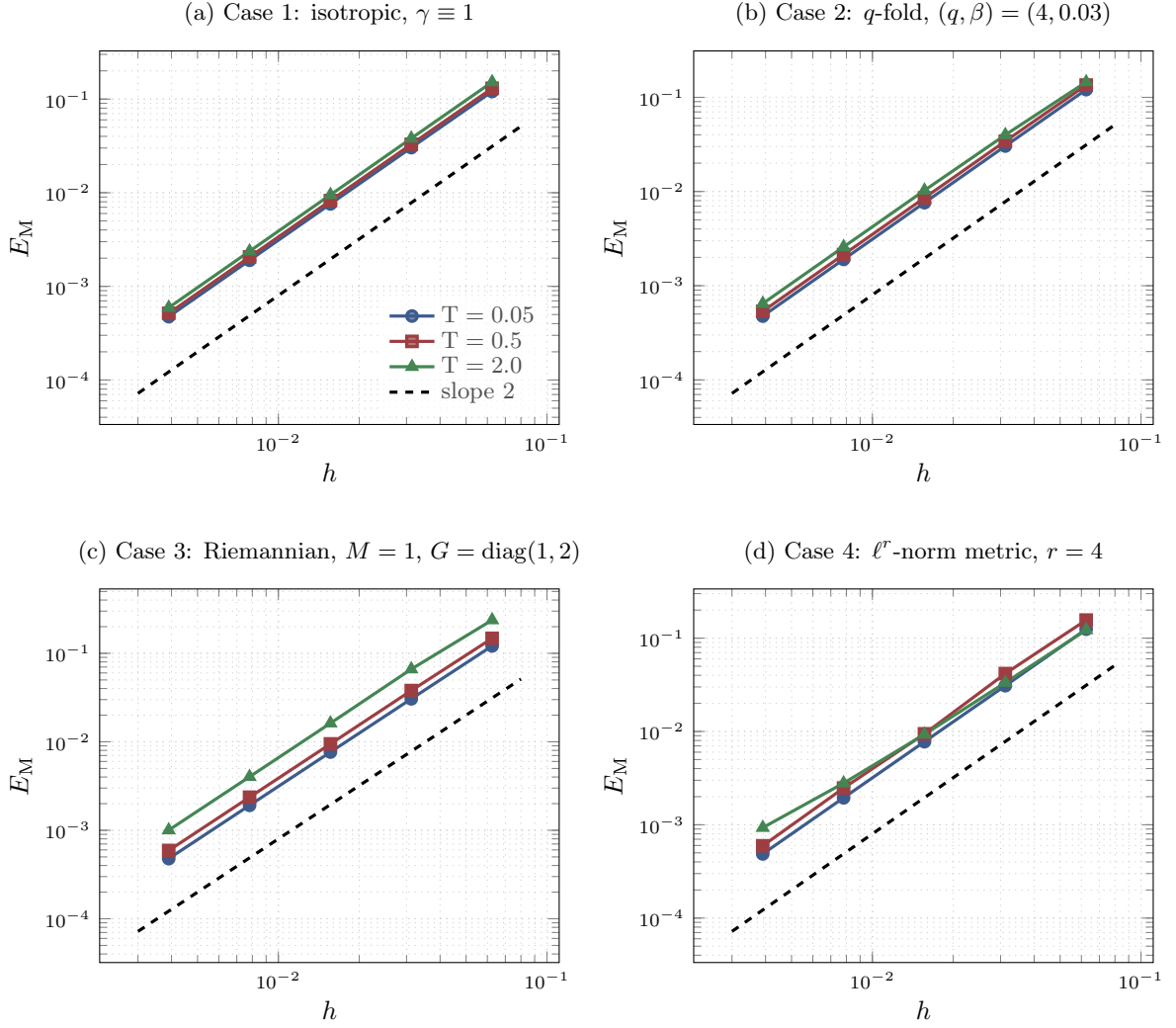
\begin{figure}[h!]
    \centering
    \begin{tikzpicture}
        \pgfplotsset{
            every axis/.append style={
                width=0.47\linewidth, height=0.40\linewidth,
                xmode=log, ymode=log,
                xlabel style={font=\normalsize, yshift=2pt},
                ylabel style={font=\normalsize, yshift=-2pt},
                title style={font=\small},
                tick label style={font=\scriptsize},
                grid=both,
                grid style={dotted, gray!40},
                every axis plot/.append style={very thick},
                legend style={font=\small, draw=none, fill opacity=0.7,
                              inner sep=2pt, row sep=-2pt},
                legend cell align=left,
            },
        }
        \begin{groupplot}[
            group style={group size=2 by 2,
                         horizontal sep=1.8cm, vertical sep=2.2cm},
        ]
        % --- (a) Case 1: isotropic ---
        \nextgroupplot[
            title={(a) Case 1: isotropic, $\gamma \equiv 1$},
            xlabel={$h$}, ylabel={$E_{\mathrm{M}}$},
            legend pos=south east,
        ]
            \addplot[seabornblue!80!black, mark=*]
                table[x index=0, y index=1] {figure/eoc_apcsf_iso_ellipse.dat};
                \addlegendentry{$\mathrm{T}=0.05$}
            \addplot[seabornred!80!black, mark=square*]
                table[x index=0, y index=2] {figure/eoc_apcsf_iso_ellipse.dat};
                \addlegendentry{$\mathrm{T}=0.5$}
            \addplot[seaborngreen!80!black, mark=triangle*]
                table[x index=0, y index=3] {figure/eoc_apcsf_iso_ellipse.dat};
                \addlegendentry{$\mathrm{T}=2.0$}
            \addplot[domain=0.003:0.08, dashed, black] {8*x^2};
                \addlegendentry{slope 2}
        % --- (b) Case 2: q-fold ---
        \nextgroupplot[
            title={(b) Case 2: $q$-fold, $(q,\beta)=(4,0.03)$},
            xlabel={$h$}, ylabel={$E_{\mathrm{M}}$},
            legend pos=south east,
        ]
            \addplot[seabornblue!80!black, mark=*]
                table[x index=0, y index=1] {figure/eoc_apcsf_kfold_k4_weak_ellipse.dat};
        %        \addlegendentry{$\mathrm{T}=0.05$}
            \addplot[seabornred!80!black, mark=square*]
                table[x index=0, y index=2] {figure/eoc_apcsf_kfold_k4_weak_ellipse.dat};
        %        \addlegendentry{$\mathrm{T}=0.5$}
            \addplot[seaborngreen!80!black, mark=triangle*]
                table[x index=0, y index=3] {figure/eoc_apcsf_kfold_k4_weak_ellipse.dat};
        %        \addlegendentry{$\mathrm{T}=2.0$}
            \addplot[domain=0.003:0.08, dashed, black] {8*x^2};
        %        \addlegendentry{slope 2}
        % --- (c) Case 3: Riemannian ---
        \nextgroupplot[
            title={(c) Case 3: Riemannian, $M=1$, $G=\mathrm{diag}(1,2)$},
            xlabel={$h$}, ylabel={$E_{\mathrm{M}}$},
            legend pos=south east,
        ]
            \addplot[seabornblue!80!black, mark=*]
                table[x index=0, y index=1] {figure/eoc_apcsf_riemannian_ellipse.dat};
         %       \addlegendentry{$\mathrm{T}=0.05$}
            \addplot[seabornred!80!black, mark=square*]
                table[x index=0, y index=2] {figure/eoc_apcsf_riemannian_ellipse.dat};
        %        \addlegendentry{$\mathrm{T}=0.5$}
            \addplot[seaborngreen!80!black, mark=triangle*]
                table[x index=0, y index=3] {figure/eoc_apcsf_riemannian_ellipse.dat};
        %        \addlegendentry{$\mathrm{T}=2.0$}
            \addplot[domain=0.003:0.08, dashed, black] {8*x^2};
        %        \addlegendentry{slope 2}
        % --- (d) Case 4: l^r-norm metric ---
        \nextgroupplot[
            title={(d) Case 4: $\ell^r$-norm metric, $r=4$},
            xlabel={$h$}, ylabel={$E_{\mathrm{M}}$},
            legend pos=south east,
        ]
            \addplot[seabornblue!80!black, mark=*]
                table[x index=0, y index=1] {figure/eoc_apcsf_l4_ellipse.dat};
        %        \addlegendentry{$\mathrm{T}=0.05$}
            \addplot[seabornred!80!black, mark=square*]
                table[x index=0, y index=2] {figure/eoc_apcsf_l4_ellipse.dat};
        %        \addlegendentry{$\mathrm{T}=0.5$}
            \addplot[seaborngreen!80!black, mark=triangle*]
                table[x index=0, y index=3] {figure/eoc_apcsf_l4_ellipse.dat};
        %        \addlegendentry{$\mathrm{T}=2.0$}
            \addplot[domain=0.003:0.08, dashed, black] {8*x^2};
        %       \addlegendentry{slope 2}
        \end{groupplot}
    \end{tikzpicture}
    \caption{Adjacent-grid manifold distance
    $E_{\mathrm{M}}$ for the AP-CSF. The initial curve is the
    $2{:}1$ ellipse. The results are shown at $\mathrm{T}=0.05$, $0.5$,
    and $2.0$ for the isotropic density and three weakly
    anisotropic densities. The dashed lines have slope two.}
    \label{fig:eoc_apcsf}
\end{figure}

\iffalse
\begin{figure}[h!]
    \centering
    \begin{subfigure}[t]{0.245\textwidth}
        \centering
        \includegraphics[width=\linewidth]{figure/wulff_compare_iso.pdf}
        \caption{Case 1}
        \label{fig:wulff_a}
    \end{subfigure}\hfill%
    \begin{subfigure}[t]{0.245\textwidth}
        \centering
        \includegraphics[width=\linewidth]{figure/wulff_compare_kfold.pdf}
        \caption{Case 2}
        \label{fig:wulff_b}
    \end{subfigure}\hfill%
    \begin{subfigure}[t]{0.245\textwidth}
        \centering
        \includegraphics[width=\linewidth]{figure/wulff_compare_riemannian.pdf}
        \caption{Case 3}
        \label{fig:wulff_c}
    \end{subfigure}\hfill%
    \begin{subfigure}[t]{0.245\textwidth}
        \centering
        \includegraphics[width=\linewidth]{figure/wulff_compare_l4.pdf}
        \caption{Case 4}
        \label{fig:wulff_d}
    \end{subfigure}
    \caption{Comparison of the long-time AP-CSF profiles computed
    by the MMS--PFEM (dashed red) with the corresponding theoretical
    Wulff profiles (solid black). The four surface energy densities
    are those used in Fig. \ref{fig:eoc_apcsf}.}
    \label{fig:wulff}
\end{figure}
\fi

The CSF results in Fig. \ref{fig:eoc_csf} exhibit approximately
second-order decay over the displayed mesh range. This behavior is
observed for the isotropic density and all three weakly anisotropic
densities. Under the coupling $\tau=h^2$, these results are consistent
with an $O(h^2+\tau)$ error for smooth curve evolutions. The AP-CSF results in Fig. \ref{fig:eoc_apcsf} show a similar trend.
The observed slopes remain close to two over the displayed mesh
range. Thus, the scalar Lagrange multiplier does not reduce the
observed coupled convergence order in these tests.

%Fig.~\ref{fig:wulff} provides an additional qualitative test. The
%long-time AP-CSF profiles agree closely with the corresponding Wulff
%profiles after the normalization used for the comparison. This
%behavior is consistent with the expected equilibrium shapes of the
%area-constrained anisotropic energies.

\subsection{Oscillations in the strongly anisotropic case}
\label{sec:strong-aniso-conv}

We next consider the strongly anisotropic $q$-fold density
\[\gamma(N) = 1 + \beta\cos(q\theta),\qquad q=2, \qquad \beta = 0.4.\]
For these parameters,
$\beta(q^2 - 1) = 1.2 > 1$. We therefore use the convex splitting
\[\gamma=\gamma_c - \gamma_e,\qquad \gamma_e=C, \qquad \gamma_c=\gamma + C,\]
as described in \Cref{sec:mms-strong}.
By \Cref{rem:splitting}, the smallest constant that makes
$\gamma_c$ weakly anisotropic is
\[
C_{\min}
=\max\bigl\{0,\beta(q^2-1)-1\bigr\}
=0.2.
\]

\begin{figure}[htbp]
     \centering
     \begin{subfigure}[b]{0.75\textwidth}
         \centering
         \includegraphics[width=1.0\linewidth]{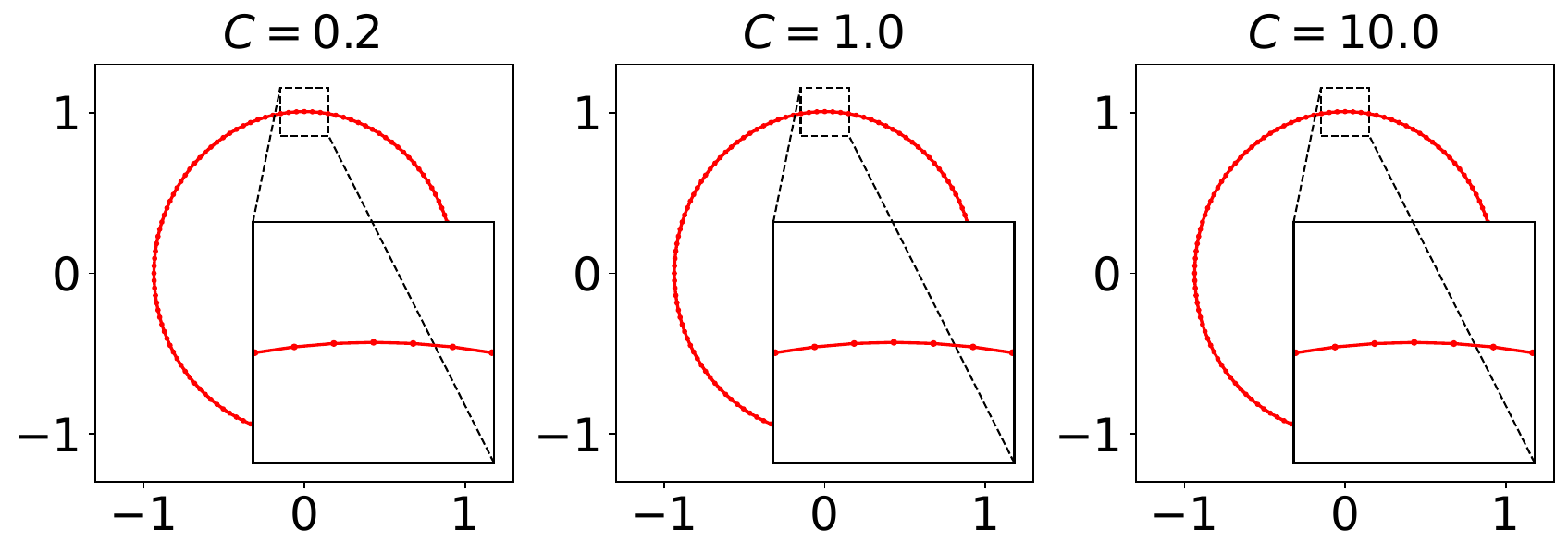}
         \caption{Mesh points $J = 128$ (no visible oscillation).}
     \end{subfigure}
     \begin{subfigure}[b]{0.75\textwidth}
         \centering
         \includegraphics[width=1\linewidth]{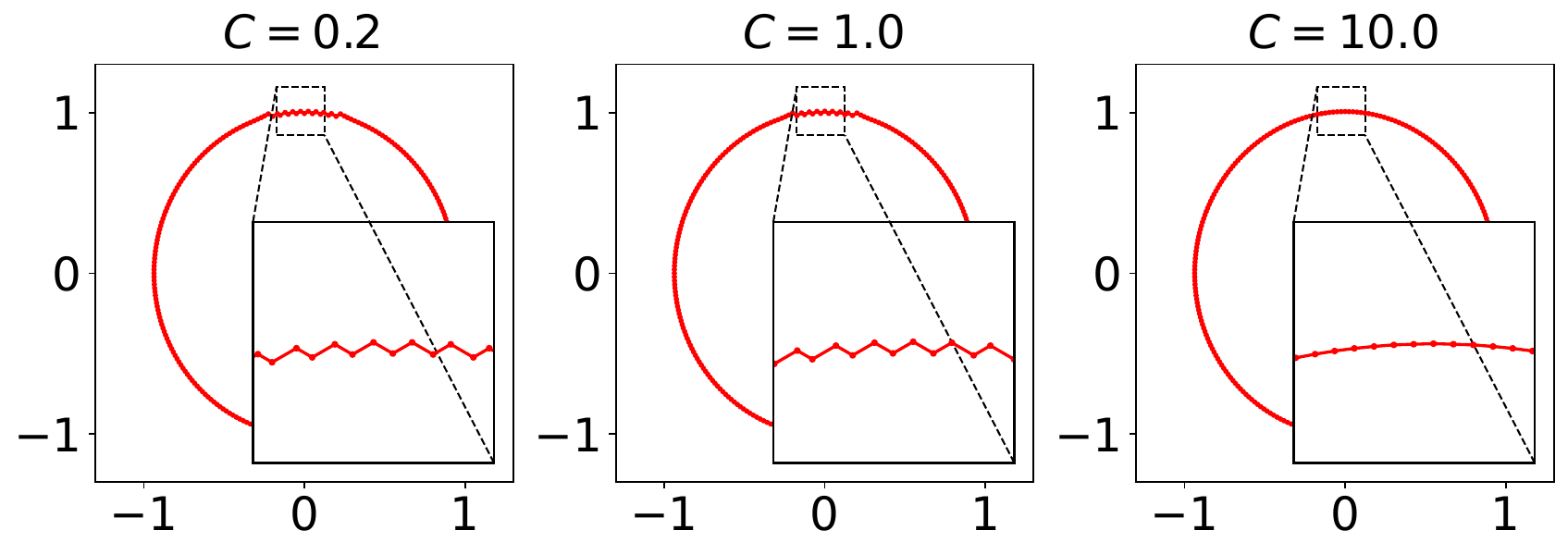}
         \caption{Mesh points $J = 256$ (visible oscillation appears for $C=0.2, 1.0$).}
     \end{subfigure}
     \begin{subfigure}[b]{0.75\textwidth}
         \centering
         \includegraphics[width=1\linewidth]{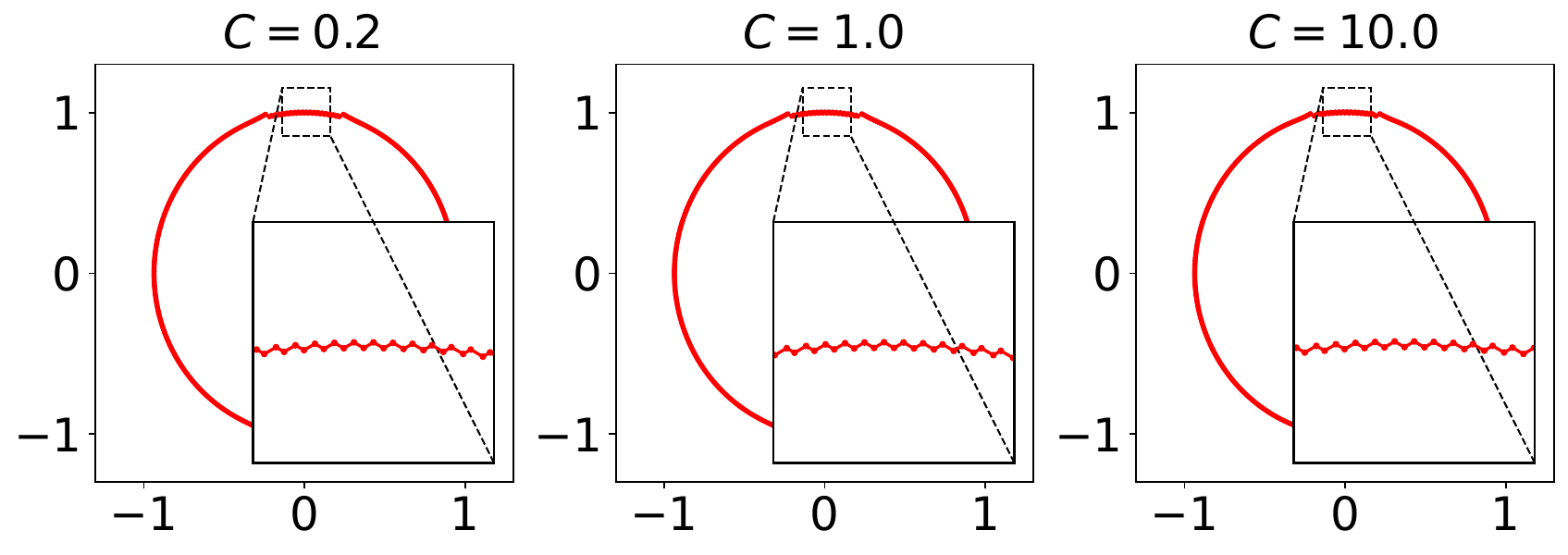}
         \caption{Mesh points $J = 512$ (oscillation appears for all  values of $C$).}
     \end{subfigure}
      \caption{Numerical curves at $\mathrm{T}=0.03$ for the strongly
    anisotropic flow with $(q,\beta)=(2,0.4)$. The three columns
    correspond, from left to right, to $C=0.2$, $1.0$, and $10.0$.}
      \label{fig:cusp-comparison}
\end{figure}

We compare three splitting constants,
$C=0.2$, $C=1.0$, and $C=10.0$.
All experiments start from the unit circle and use $\tau=h^2$.
Table~\ref{tab:eoc-strong} reports the adjacent-grid
errors at the terminal time $\mathrm{T}=0.03$, and
Fig.~\ref{fig:cusp-comparison} compares the corresponding curves at the same time.

\begin{table}[h!]
\centering
\caption{Adjacent-grid manifold distance and observed orders at terminal
time $\mathrm{T}=0.03$ for the strongly anisotropic $2$-fold density
with $(q,\beta)=(2,0.4)$.}
\label{tab:eoc-strong}
\renewcommand{\arraystretch}{1.15}
\begin{tabular}{c|cc|cc|cc}
\hline
& \multicolumn{2}{c|}{$C=0.2$}
& \multicolumn{2}{c|}{$C=1.0$}
& \multicolumn{2}{c}{$C=10.0$} \\
$(N_c,N_f)$
& $E_{\mathrm{M}}$ & Order
& $E_{\mathrm{M}}$ & Order
& $E_{\mathrm{M}}$ & Order \\
\hline
$(16,32)$
& 5.907e-02 & ---
& 5.860e-02 & ---
& 5.492e-02 & --- \\
$(32,64)$
& 1.484e-02 & 1.993
& 1.472e-02 & 1.993
& 1.377e-02 & 1.996 \\
$(64,128)$
& 3.708e-03 & 2.001
& 3.676e-03 & 2.002
& 3.435e-03 & 2.003 \\
$(128,256)$
& 7.567e-03 & $-$1.029
& 6.579e-03 & $-$0.840
& 8.709e-04 & 1.980 \\
$(256,512)$
& 6.796e-03 & 0.155
& 7.907e-03 & $-$0.265
& 5.440e-03 & $-$2.643 \\
\hline
\end{tabular}
\end{table}

At $\mathrm{T}=0.03$,
Table~\ref{tab:eoc-strong} reveals two distinct mesh regimes.
On the coarser mesh pairs, the errors decrease at an observed order
close to two. On finer mesh pairs, the errors cease to decrease
monotonically. Several orders then become negative. Consequently, the
complete refinement sequences do not demonstrate mesh convergence.
They instead show a coarse-grid slope-two regime followed by a loss
of adjacent-grid agreement.
The splitting constant affects the mesh level at
which this loss becomes visible: for $C=0.2$ and $C=1.0$, the
deterioration begins at the mesh pair $(128,256)$, whereas for
$C=10.0$ it is delayed until $(256,512)$. Thus, increasing $C$ can
delay the deterioration, but it does not eliminate it.

Fig.~\ref{fig:cusp-comparison} associates the loss of adjacent-grid
agreement with increasingly sharp, mesh-scale oscillations.
The present data do not identify a unique mechanism for this
behavior. Possible causes include the loss of parabolicity of the
continuum law, the splitting-dependent discrete dynamics, unresolved
small scales, or a combination of these effects. Therefore, this
experiment compares finite-dimensional discrete solutions rather
than establishing convergence to a strongly anisotropic continuum
flow.
The parameter $C$ should not be interpreted as a physical isotropic
shortening coefficient. The implicit and explicit contributions of
the splitting cancel at the level of the original anisotropic energy,
as discussed after \Cref{thm:mms-strong}. For every
$C\geq C_{\min}$, the one-step energy estimate remains valid for a
global minimizer of the discrete problem.
%A converged KKT residual alone does not establish that the computed stationary point is such a minimizer.

\subsection{Numerical assessment of the structure-preserving properties}
\label{sec:acsf-evolution}

We next examine the discrete energy decay and mesh behavior of the anisotropic CSF schemes. Unless stated otherwise, the initial curve
is the unit circle. We monitor the normalized anisotropic energy and the mesh ratio
\[
\l.\frac{\cE_{\gamma}(t)}{\cE_{\gamma}(0)}\r|_{t=t_k}
:=
\frac{\cE_\gamma(X_h^k)}{\cE_\gamma(X_h^0)},\quad \l.\Psi(t)\r|_{t=t_k}
:=\frac{\displaystyle\max_{1\leq j\leq J}|\sigma_j^k|}
     {\displaystyle\min_{1\leq j\leq J}|\sigma_j^k|},
\qquad
\Gamma_h^k=\bigcup_{j=1}^J\sigma_j^k,
\]
where $|\sigma_j^k|$ is the length of the $j$th polygonal edge.
The quantity $\Psi$ measures the nonuniformity of the mesh and
satisfies $\Psi(t)\geq1$. We also record the number of Newton iterations at each time step. In the following experiments, we use
$J=128$, $\tau=5\times10^{-4}$.

\begin{figure}[h!]
     \begin{subfigure}[b]{\textwidth}
         \centering
         \includegraphics[width=0.8\linewidth,height=0.7\linewidth]{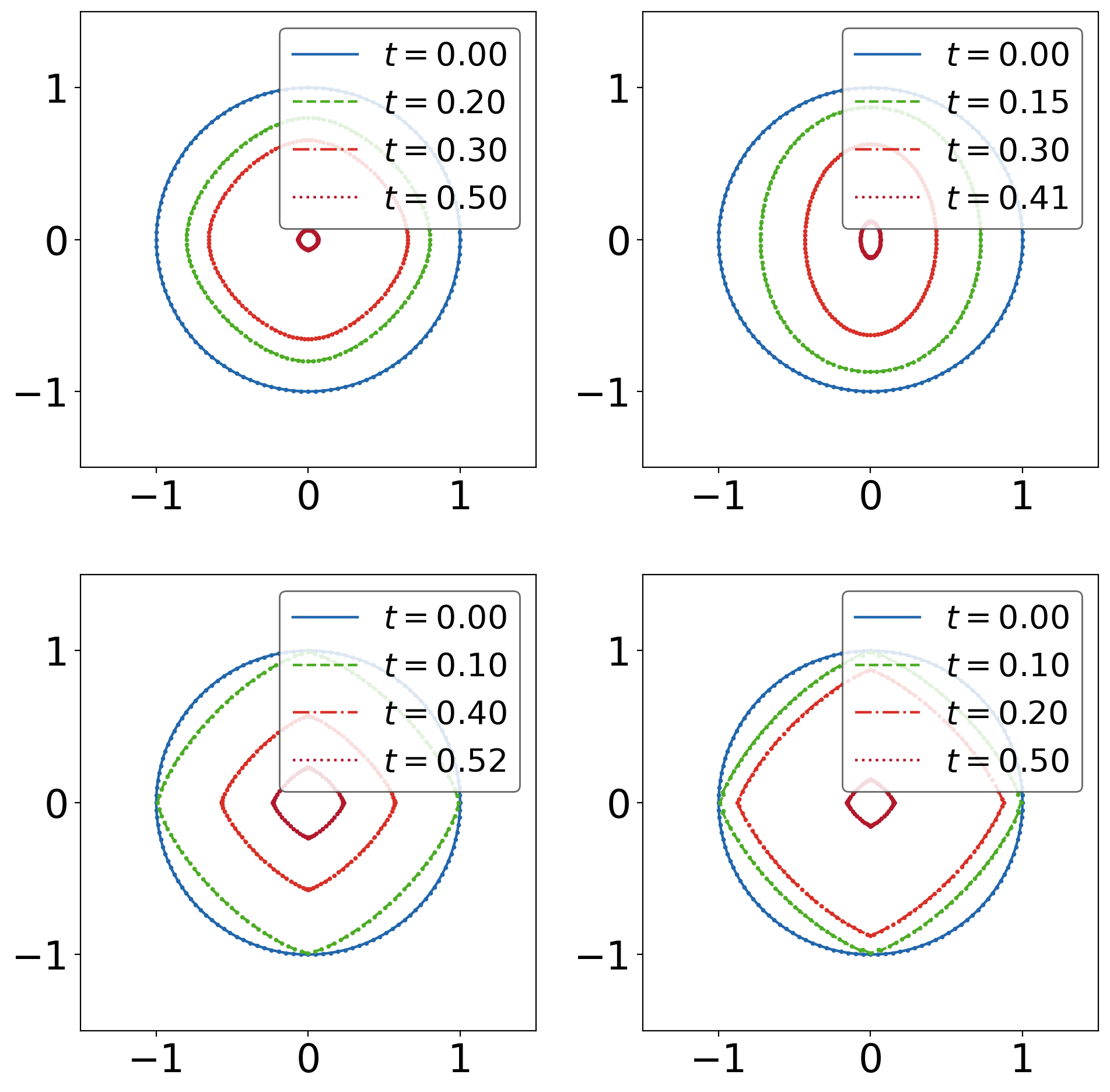}
         \caption{Snapshots under four surface energy densities.
       Top left: $q$-fold density with
        $(q,\beta)=(4,0.03)$.
        Top right: Riemannian density with
        $M=1$ and $G=\mathrm{diag}(1,2)$.
        Bottom left: $\ell^4$-norm density.
        Bottom right: strongly anisotropic $q$-fold density with
        $(q,\beta)=(4,0.1)$ and $C=0.5$.}
     \end{subfigure}
     \begin{subfigure}[b]{\textwidth}
         \centering
         \includegraphics[width=0.95\linewidth,height=0.35\linewidth]{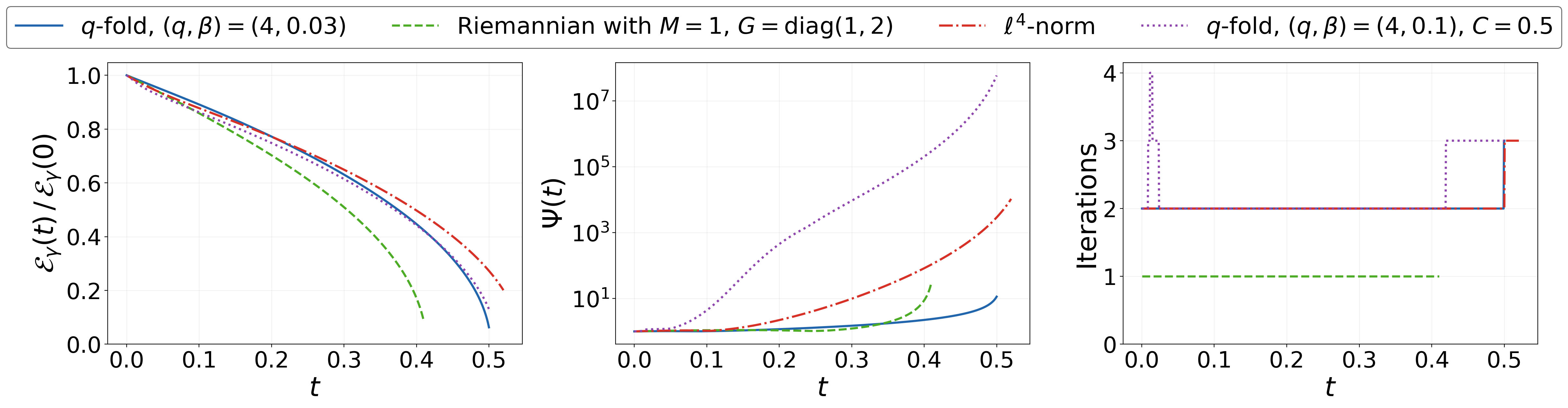}
         \caption{Normalized energy, mesh ratio, and number of Newton
        iterations.}
     \end{subfigure}
     \vspace{-5mm}
      \caption{Numerical evolution of the unit circle under four anisotropic CSFs.}
      \label{fig:csf_evolution}
\end{figure}

Fig.~\ref{fig:csf_evolution} compares three weakly anisotropic
densities from \Cref{Example:Ani} and one strongly anisotropic
density. The weakly anisotropic examples are the $q$-fold density
with $(q,\beta)=(4,0.03)$, the Riemannian density with $G=\mathrm{diag}(1,2)$, and the $\ell^4$-norm density
$\gamma(N)=(n_1^4+n_2^4)^{1/4}$.
The strongly anisotropic example uses the $q$-fold density with
$(q,\beta)=(4,0.1)$. For this case, we take
$C=C_{\min}=0.5$ in the convex-splitting scheme.

In the weakly anisotropic examples, the numerical curves become
increasingly close in shape to the corresponding Wulff profiles while
they shrink. The normalized energy decreases monotonically, in
agreement with the one-step energy estimate. The mesh ratios remain
moderate over the displayed time intervals. Newton's method requires
only a few iterations per time step.
In the strongly anisotropic example, sharp facet-like features
develop in directions determined by the anisotropy. The recorded
energy still decreases monotonically, which is consistent with
\eqref{eq:energy-aniso-strong}. However, the mesh ratio increases
markedly as the mesh becomes nonuniform near these sharp features.
Newton's method nevertheless satisfies the prescribed residual
tolerance at every reported step.

\begin{figure}[h!]
     \begin{subfigure}[b]{\textwidth}
         \centering
         \includegraphics[width=1.0\linewidth]{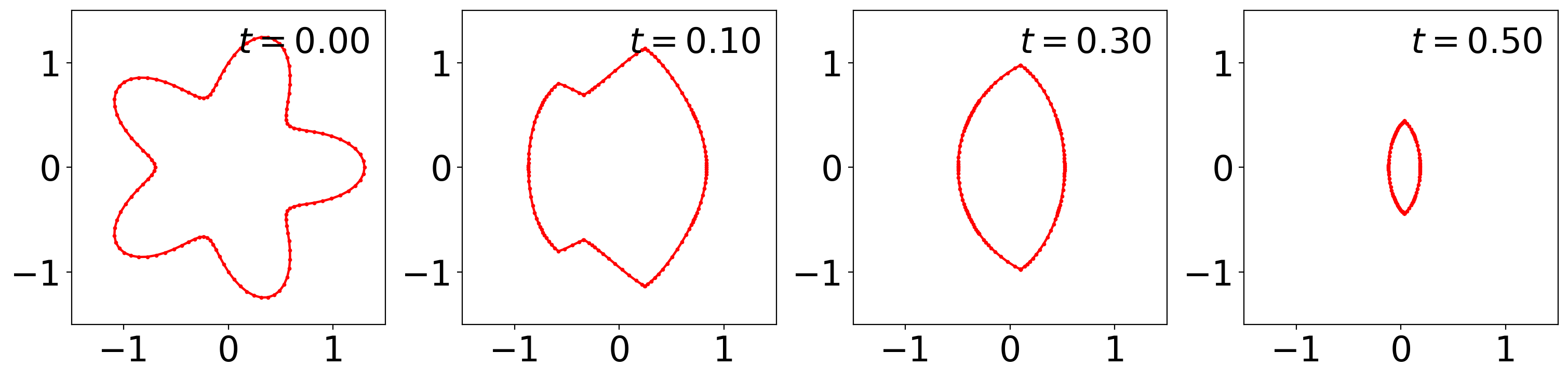}
         \caption{Snapshots of the curve evolution.}
     \end{subfigure}
     \begin{subfigure}[b]{\textwidth}
         \centering
         \includegraphics[width=1\linewidth]{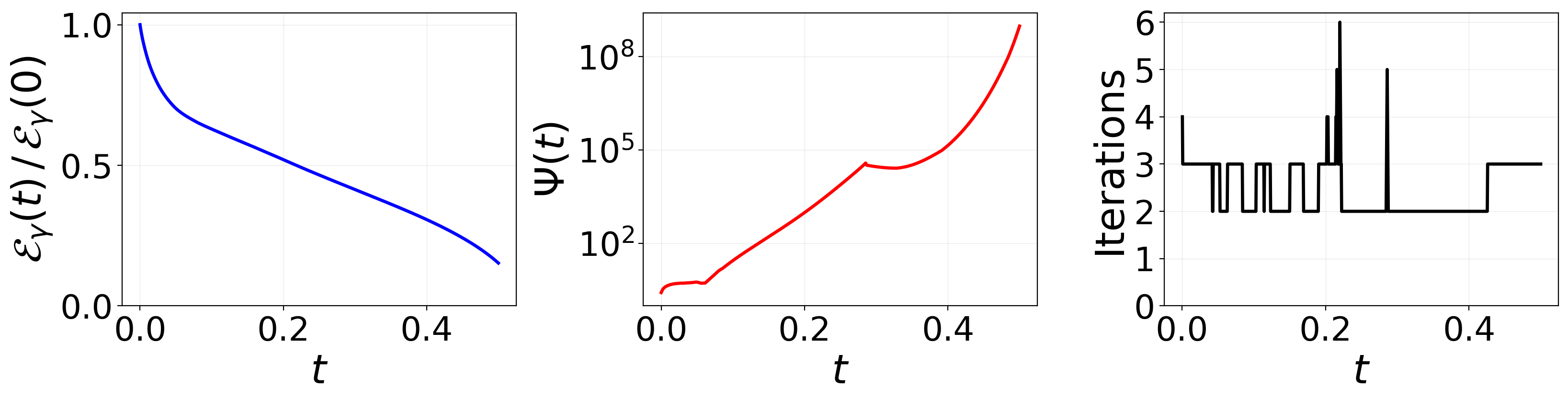}
         \caption{Normalized energy, mesh ratio, and number of Newton iterations.}
     \end{subfigure}
      \caption{Evolution of a ``flower" curve under the strongly
    anisotropic $q$-fold density with
    $(q,\beta)=(2,0.4)$ and $C=0.2$.}
      \label{fig:csf_case2_q2_strong_flower}
\end{figure}

We also consider a nonconvex initial curve. The initial shape is a
``flower'' curve, and the anisotropy is the strongly anisotropic $q$-fold density with
$(q,\beta)=(2,0.4)$.
As shown in Fig. \ref{fig:csf_case2_q2_strong_flower}, the petals
initially contract and sharp facets subsequently emerge. The curve
then continues to shrink. The recorded discrete energy decreases
throughout the computation, and Newton's method converges in a small
number of iterations at each time step.

\begin{figure}[h!]
     \begin{subfigure}[b]{\textwidth}
         \centering
         \includegraphics[width=1.0\linewidth]{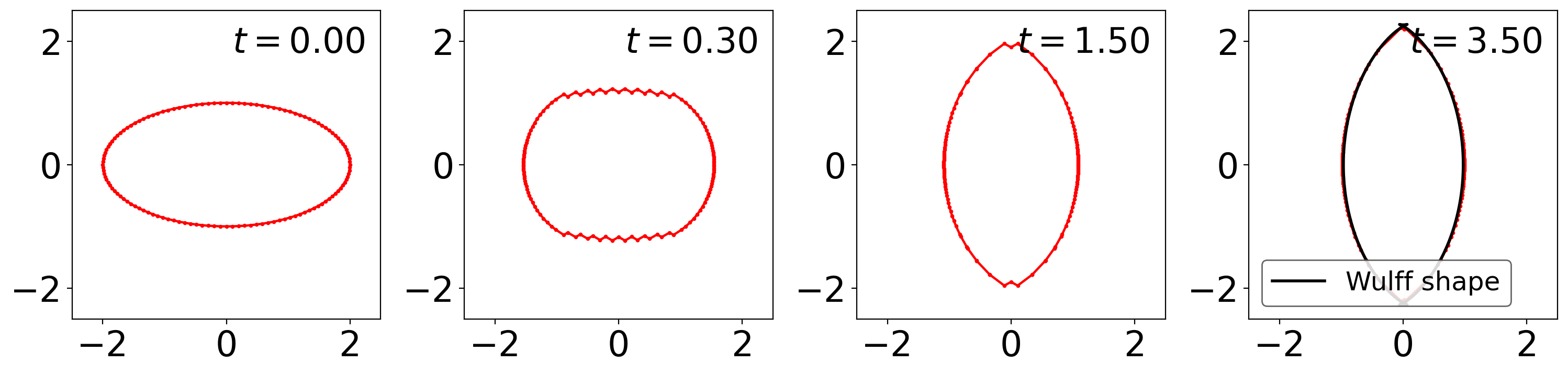}
         \caption{Snapshots of the curve evolution.}
     \end{subfigure}
     \begin{subfigure}[b]{\textwidth}
         \centering
         \includegraphics[width=1\linewidth]{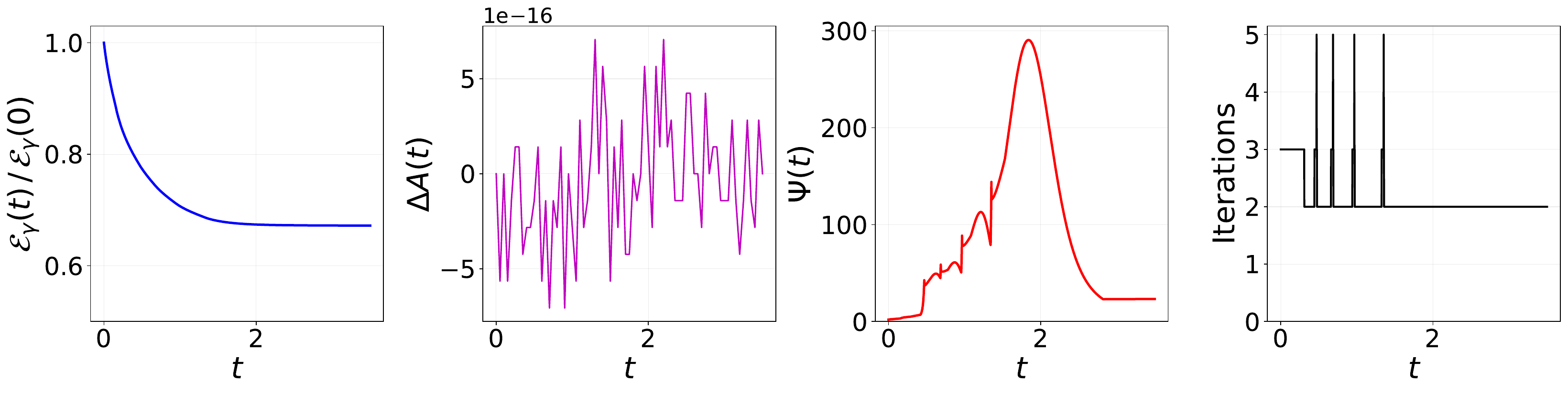}
         \caption{Normalized energy, relative area loss, mesh ratio, and number of Newton iterations.}
     \end{subfigure}
      \caption{Evolution of an ellipse under the strongly anisotropic AP-CSF with $(q,\beta)=(2,0.4)$ and $C=0.2$.}
      \label{fig:apcsf_q2_strong}
\end{figure}

\section{Conclusion}\label{sec:conclusion}

In this paper, we developed a variational MMS--PFEM framework for isotropic and anisotropic curve shortening flows. The framework combines a minimizing-movement formulation compatible with parametric finite elements and admissible modified energies. For any positive surface energy density $\gamma\in C^2(\mathbb S^1)$, the unconstrained scheme admits a unique solution at each time step for an arbitrary time-step size, provided that the parametric derivative of the current discrete curve does not vanish on any mesh element.
Furthermore, we rigorously prove that the proposed schemes satisfy the energy dissipation property,
no matter whether the surface energy density is isotropic, weakly anisotropic or strongly anisotropic.
The framework recovers Dziuk's isotropic scheme, applies directly to weak anisotropies, and treats strong anisotropies through convex splitting. A single scalar constraint further extends the method to area-preserving flows. Each constrained step admits a global minimizer that preserves the signed area exactly and dissipates the energy, although uniqueness is not guaranteed. The numerical experiments exhibit the expected convergence for smooth isotropic and weakly anisotropic examples, faithfully capture the facet-like profiles but may develop mesh-scale oscillations and cusps under refinement for the strongly anisotropic cases. They also confirm favorable mesh quality, energy decay and area preservation.
 Future work will extend the framework to surface diffusion and geometric surface flows in $\mathbb{R}^3$.

\section*{Acknowledgement}
This work is supported by the NSFC 12671530 and 12271414 (Wei Jiang), NSFC 12201342 (Chunmei Su) and a fellowship award from the Research Grants Council of the Hong Kong Special Administrative Region PolyU/RFS2324-5S03 (Ganghui Zhang).

%\section*{Declaration of competing interest}
%The authors declare that they have no competing interests.

\bibliographystyle{elsarticle-num}

\bibliography{ref}

\end{document}